\documentclass[12pt]{article}%

\usepackage{amsmath}
\usepackage{amssymb}
\usepackage{epsfig}
\usepackage{amscd}
\usepackage{amsfonts}
\usepackage{graphicx}%

\newtheorem{theorem}{Theorem}

\newtheorem{corollary}{Corollary}

\newtheorem{proposition}{Proposition}
\newtheorem{remark}{Remark}
\numberwithin{equation}{section}
\numberwithin{lemma}{section}
\numberwithin{theorem}{section}
\numberwithin{remark}{section}
\numberwithin{corollary}{section}
\numberwithin{proposition}{section}
\numberwithin{definition}{section}
\numberwithin{table}{section}

\begin{document}

\title{Layer methods for Navier-Stokes equations with additive noise}

\author{G.N. Milstein\thanks{Ural State University, Lenin Str.~51, 620083
Ekaterinburg, Russia; email: Grigori.Milstein@usu.ru} 
\and M.V. Tretyakov\thanks{School of Mathematical Sciences, University of
Nottingham, University Park, Nottingham, NG7 2RD, UK, email: Michael.Tretyakov@nottingham.ac.uk}}
\maketitle

\begin{abstract}
We propose and study a number of layer methods for stochastic Navier-Stokes
equations (SNSE) with spatial periodic boundary conditions and additive noise.
The methods are constructed using conditional probabilistic representations of
solutions to SNSE and exploiting ideas of the weak sense numerical integration
of stochastic differential equations. We prove some convergence results for
the proposed methods. Results of numerical experiments on two model problems
are presented.



\end{abstract}

\section{Introduction}

Let $(\Omega,\mathcal{F},P)$ be a probability space and $(w(t),\mathcal{F}%
_{t}^{w})=((w_{1}(t),\ldots,w_{q}(t))^{\top},\mathcal{F}_{t}^{w})$ be a
$q$-dimensional standard Wiener process, where $\mathcal{F}_{t}^{w},\ 0\leq
t\leq T,$ is an increasing family of $\sigma$-subalgebras of $\mathcal{F}$
induced by $w(t).$ We consider the system of stochastic Navier-Stokes
equations (SNSE) with additive noise for velocity $v$ and pressure $p$ in a
viscous incompressible flow:%
\begin{gather}
dv(t)=\left[  \frac{\sigma^{2}}{2}\Delta v-(v,\nabla)v-\nabla p+f(t,x)\right]
dt+\sum_{r=1}^{q}\gamma_{r}(t,x)dw_{r}(t),\label{NS1}\\
\ \ 0\leq t\leq T,\ x\in \mathbf{R}^{n}, \nonumber \\
\operatorname{div}\ v=0, \label{NS2}%
\end{gather}
with spatial periodic conditions%
\begin{align}
v(t,x+Le_{i})  &  =v(t,x),\ p(t,x+Le_{i})=p(t,x),\label{NS4}\\
0  &  \leq t\leq T,\ \ i=1,\ldots,n,\nonumber
\end{align}
and the initial condition%
\begin{equation}
v(0,x)=\varphi(x). \label{NS3}%
\end{equation}
In (\ref{NS1})-(\ref{NS2}) we have $v\in\mathbf{R}^{n}$,$\ p$ is a
scalar,\ $f$ $\in\mathbf{R}^{n},$\ $\gamma_{r}\in\mathbf{R}^{n}$; $\{e_{i}\}$
is the canonical basis in $\mathbf{R}^{n}$ and $L>0$ is the period (for
simplicity in writing, the periods in all the directions are taken the same).
The functions $f=f(t,x)$ and $\gamma_{r}(t,x)$ are supposed to be spatial
periodic as well. Further, we require that $\gamma_{r}(t,x)$ are divergence
free:
\begin{equation}
\operatorname{div}\gamma_{r}(t,x)=0,\ r=1,\ldots,q. \label{NS03}%
\end{equation}

SNSE can be useful for explaining the turbulence phenomenon (see
\cite{turb1,Flandoli,RozNS04} and references therein). They have complicated
dynamics and some interesting properties (e.g., ergodicity of solutions
\cite{HaMa06,HasNS,DaPr,MaPa06}). At the same time, rather little has been
done in numerics for SNSE. Let us cite \cite{HRoz07}, where algorithms based
on Wiener Chaos expansion are considered, and quite recent works
\cite{BCP10,CP,D}, where splitting schemes with finite element or Galerkin
approximations are applied. Here we suggest to exploit some probabilistic
representations of solutions to SNSE for constructing numerical methods of the
layer type. The proposed methods are promised to be effective, reliable
numerical methods for studying SNSE. Layer methods for deterministic
semilinear and quasilinear partial differential equations of parabolic type
were proposed in \cite{M1,quasic} (see also \cite{MT1,french}), and for
deterministic NSEs they were first considered in \cite{BM} and further
developed in \cite{NS5,NSB}. Layer methods for linear and semilinear
stochastic partial differential equations (SPDE) were constructed and analyzed
in \cite{spde}.

The rest of the paper is organized as follows. In Section~\ref{Secpre} we
introduce additional notation and write down probabilistic representations for
linearized SNSE (i.e., stochastic Oseen-Stokes equations) and for the SNSE
(\ref{NS1})-(\ref{NS3}) which we use in Section~\ref{secLayer} for
constructing layer methods for the SNSE. Three layer methods are given in
Section~\ref{secLayer} together with discussion of their implementation.
Numerical error analysis is done in Section~\ref{secER}. Results of numerical
experiments on two test models are presented in Section~\ref{secnum}.

\section{Preliminaries\label{Secpre}}

In this section we recall the required function spaces
\cite{CM,RT,T,RozNS04,RozNS05} and write probabilistic representations of
solutions to linearized SNSE and to SNSE resting on results from
\cite{KryR,Kun,Pard,R}.

\subsection{Function spaces, the Helmholtz-Hodge-Leray decomposition, and
notation}

Let $\{e_{i}\}$ be the canonical basis in $\mathbf{R}^{n}.$ We shall consider
spatial periodic $n$-vector functions $u(x)=(u^{1}(x),\ldots,u^{n}(x))^{\top}$
in $\mathbf{R}^{n}:$ $u(x+Le_{i})=u(x),\ i=1,\ldots,n,$ where $L>0$ is the
period in $i$th direction. Denote by $Q=(0,L)^{n}$ the cube of the period (of
course, one may consider different periods $L_{1},\ldots,L_{n}$ in the
different directions).\ We denote by $\mathbf{L}^{2}(Q)$ the Hilbert space of
functions on $Q$ with the scalar product and the norm%
\[
(u,v)=\int_{Q}\sum_{i=1}^{n}u^{i}(x)v^{i}(x)dx,\ \Vert u\Vert=(u,u)^{1/2}.
\]
We keep the notation $|\cdot|$ for the absolute value of numbers and for the
length of $n$-dimensional vectors, for example,%
\[
|u(x)|=[(u^{1}(x))^{2}+\cdots+(u^{n}(x))^{2}]^{1/2}.
\]

We denote by $\mathbf{H}_{p}^{m}(Q),\ m=0,1,\ldots,$ the Sobolev space of
functions which are in $\mathbf{L}^{2}(Q),$ together with all their
derivatives of order less than or equal to $m,$ and which are periodic
functions with the period $Q.$ The space $\mathbf{H}_{p}^{m}(Q)$ is a Hilbert
space with the scalar product and the norm%
\[
(u,v)_{m}=\int_{Q}\sum_{i=1}^{n}\sum_{[\alpha^{i}]\leq m}D^{\alpha^{i}}%
u^{i}(x)D^{\alpha^{i}}v^{i}(x)dx,\ \Vert u\Vert_{m}=[(u,u)_{m}]^{1/2},
\]
where $\alpha^{i}=(\alpha_{1}^{i},\ldots,\alpha_{n}^{i}),\ \alpha_{j}^{i}%
\in\{0,\ldots,m\},\ [\alpha^{i}]=\alpha_{1}^{i}+\cdots+\alpha_{n}^{i},$ and%
\[
D^{\alpha^{i}}=D_{1}^{\alpha_{1}^{i}}\cdots D_{n}^{\alpha_{n}^{i}}%
=\frac{\partial^{\lbrack\alpha^{i}]}}{(\partial x^{1})^{\alpha_{1}^{i}}%
\cdots(\partial x^{n})^{\alpha_{n}^{i}}}\ ,\ i=1,\ldots,n.
\]
Note that $\mathbf{H}_{p}^{0}(Q)=\mathbf{L}^{2}(Q).$

Introduce the Hilbert subspaces of $\mathbf{H}_{p}^{m}(Q):$%
\begin{align*}
\mathbf{V}_{p}^{m}  &  =\{v:\ v\in\mathbf{H}_{p}^{m}(Q),\ \operatorname{div}%
v=0\},\ m>0,\\
\mathbf{V}_{p}^{0}  &  =\text{the closure of }\mathbf{V}_{p}^{m},\ m>0\text{
in }\mathbf{L}^{2}(Q).
\end{align*}
Clearly,%
\[
\mathbf{V}_{p}^{m_{1}}=\text{the closure of }\mathbf{V}_{p}^{m_{2}}\text{ in
}\mathbf{H}_{p}^{m_{1}}(Q)\text{ for any}\ m_{2}\geq m_{1}.
\]

Denote by $P$ the orthogonal projection in $\mathbf{H}_{p}^{m}(Q)$ onto
$\mathbf{V}_{p}^{m}$ (we omit $m$ in the notation $P$ here). The operator $P$
is often called the Leray projection. Due to the Helmholtz-Hodge-Leray
decomposition, any function $u\in\mathbf{H}_{p}^{m}(Q)$ can be represented as
\[
u=Pu+\nabla g,\ \operatorname{div}Pu=0,
\]
where $g=g(x)$ is a scalar $Q$-periodic function such that $\nabla
g\in\mathbf{H}_{p}^{m}(Q).$ It is natural to introduce the notation $P^{\bot
}u:=\nabla g$ and hence write
\[
u=Pu+P^{\bot}u
\]
with
\[
P^{\bot}u\in(\mathbf{V}_{p}^{m})^{\bot}=\{v:\ v\in\mathbf{H}_{p}%
^{m}(Q),\ v=\nabla g\}.
\]

Let%
\begin{gather}
u(x)=\sum_{\mathbf{n}\in\mathbf{Z}^{n}}u_{\mathbf{n}}e^{i(2\pi/L)(\mathbf{n}%
,x)},\ g(x)=\sum_{\mathbf{n}\in\mathbf{Z}^{n}}g_{\mathbf{n}}e^{i(2\pi
/L)(\mathbf{n},x)},\ g_{\mathbf{0}}=0,\label{N00}\\
Pu(x)=\sum_{\mathbf{n}\in\mathbf{Z}^{n}}(Pu)_{\mathbf{n}}e^{i(2\pi
/L)(\mathbf{n},x)},\ P^{\bot}u(x)=\nabla g(x)=\sum_{\mathbf{n}\in
\mathbf{Z}^{n}}(P^{\bot}u)_{\mathbf{n}}e^{i(2\pi/L)(\mathbf{n},x)}\nonumber
\end{gather}
be the Fourier expansions of $u,$\ $g,$\ $Pu,$ and $P^{\bot}u=\nabla g.$ Here
$u_{\mathbf{n}},$\ $(Pu)_{\mathbf{n}},\ $and $(P^{\bot}u)_{\mathbf{n}}=(\nabla
g)_{\mathbf{n}}$ are $n$-dimensional vectors and $g_{\mathbf{n}}$ are scalars.
We note that $g_{\mathbf{0}}$ can be any real number but for definiteness we
set $g_{\mathbf{0}}=0.$ The coefficients $(Pu)_{\mathbf{n}},\ (P^{\bot
}u)_{\mathbf{n}}$, and $g_{\mathbf{n}}$ can be easily expressed in terms of
$u_{\mathbf{n}}:$
\begin{align}
(Pu)_{\mathbf{n}}  &  =u_{\mathbf{n}}-\frac{u_{\mathbf{n}}^{\top}\mathbf{n}%
}{|\mathbf{n}|^{2}}\mathbf{n,\ }(P^{\bot}u)_{\mathbf{n}}=i\frac{2\pi}%
{L}g_{\mathbf{n}}\mathbf{n=}\frac{u_{\mathbf{n}}^{\top}\mathbf{n}}%
{|\mathbf{n}|^{2}}\mathbf{n,\ }\label{N01}\\
g_{\mathbf{n}}  &  =-i\frac{L}{2\pi}\frac{u_{\mathbf{n}}^{\top}\mathbf{n}%
}{|\mathbf{n}|^{2}},\ \mathbf{n\neq0,\ }g_{\mathbf{0}}=0.\nonumber
\end{align}

We have%
\[
\nabla e^{i(2\pi/L)(\mathbf{n},x)}=\mathbf{n}e^{i(2\pi/L)(\mathbf{n},x)}\cdot
i\frac{2\pi}{L},
\]
hence $u_{\mathbf{n}}e^{i(2\pi/L)(\mathbf{n},x)}\in\mathbf{V}_{p}^{m}$ if and
only if $(u_{\mathbf{n}},\mathbf{n)}=0.$ We obtain from here that the
orthogonal basis of the subspace $(\mathbf{V}_{p}^{m})^{\bot}$ consists of
$\mathbf{n}e^{i(2\pi/L)(\mathbf{n},x)},\ \mathbf{n}\in\mathbf{Z}%
^{n},\ \mathbf{n\neq0}$; and an orthogonal basis of $\mathbf{V}_{p}^{m}$
consists of $_{k}u_{\mathbf{n}}e^{i(2\pi/L)(\mathbf{n},x)},$\ $k=1,\ldots
,n-1,\ \mathbf{n}\in\mathbf{Z}^{n},\ $where under $\mathbf{n\neq0}$ the
vectors $_{k}u_{\mathbf{n}}$ are orthogonal to $\mathbf{n:}$\textbf{\ }%
$\mathbf{(}_{k}u_{\mathbf{n}},\mathbf{n)}=0,\ k=1,\ldots,n-1,$ and they are
orthogonal among themselves: $\mathbf{(}_{k}u_{\mathbf{n}},\ _{m}%
u_{\mathbf{n}}\mathbf{)}=0,$\ $k,m=1,\ldots,n-1,$\ $m\neq k,$ and finally, for
$\mathbf{n=0,}$ the vectors $_{k}u_{\mathbf{0}},\ k=1,\ldots,n,$ are orthogonal.

In what follows we suppose that the below assumptions hold.$\medskip$

\noindent\textbf{Assumptions 2.1.} \textit{We assume that the coefficients
}$f(t,x)$ \textit{and} $\gamma_{r}(s,x),$ $r=1,\ldots,q,$ \textit{are
sufficiently smooth} \textit{and} \textit{the problem }(\ref{NS1}%
)-(\ref{NS3})\textit{ has a unique classical solution} $v(t,x),\ p(t,x),$
$(t,x)\in\lbrack0,T]\times\mathbf{R}^{n},$ \textit{which} \textit{has
continuous derivatives in the space variable }$x$\ \textit{up to some order,
and the solution and the derivatives have uniformly in }$(t,x)$
\textit{bounded moments of a sufficiently high order }$m,$ $2\leq m<m_{0},$
\textit{where} $m_{0}>2$ \textit{is a positive number or }$m_{0}=\infty
$\textit{. }$\medskip$

The solution $v(t,x),\ p(t,x),$ $(t,x)\in\lbrack0,T]\times\mathbf{R}^{n},$ to
(\ref{NS1})-(\ref{NS3}) is $\mathcal{F}_{t}^{w}$-adaptive, $v(t,\cdot
)\in\mathbf{V}_{p}^{m}$ and $\nabla p(t,\cdot)\in(\mathbf{V}_{p}^{m})^{\bot}$
for every $t\in\lbrack0,T]$ and $\omega\in\Omega.$

Assumptions of this kind are rather usual for works dedicated to numerics.
They are rested on results concerning regularity of solutions (see, e.g., the
corresponding theory for deterministic NSE in \cite{RT,T}). Unfortunately, we
could not find explicit results on the classical solution for SNSE in
literature. At the same time, the question about existence of the unique
sufficiently regular (with respect to $x)$ solution of the SNSE (\ref{NS1}%
)-(\ref{NS3}) on a time interval $[0,T]$ is analogous to the one in the
deterministic case. Indeed, the following remark reduces this problem of
regularity for the SNSE to regularity of solutions to NSE with random
coefficients which is close to the theory of deterministic NSE treated in
\cite{RT,T}.

\begin{remark}
Let $\Gamma(t,x)=\sum_{r=1}^{q}\int_{0}^{t}\gamma_{r}(s,x)dw_{r}(s).$ Then
$V(t,x)=v(t,x)+\Gamma(t,x)$ together with $p(t,x)$ solves the following
`usual' NSE with random coefficients:
\begin{gather*}
\frac{\partial}{\partial t}V=\frac{\sigma^{2}}{2}\Delta V-(V-\Gamma
(t,x),\nabla)(V-\Gamma(t,x))-\nabla p+f(t,x)-\frac{\sigma^{2}}{2}\Delta
\Gamma(t,x),\\
0\leq t\leq T,\ x\in\mathbf{R}^{n},\\
\operatorname{div}\ V=0,
\end{gather*}
with spatial periodic conditions%
\begin{align*}
V(t,x+Le_{i})  &  =V(t,x),\ p(t,x+Le_{i})=p(t,x),\ \\
0  &  \leq t\leq T,\ \ i=1,\ldots,n,
\end{align*}
and the initial condition%
\[
V(0,x)=\varphi(x).
\]

\end{remark}

\subsection{Probabilistic representations of solutions to linearized
SNSE\label{prepOseen}}

We start with considering a linearized version of the SNSE (\ref{NS1}%
)-(\ref{NS3}), i.e., the stochastic Oseen-Stokes equations (see
\cite{MikStokes}):
\begin{gather}
dv_{a}(t)=\left[  \frac{\sigma^{2}}{2}\Delta v_{a}-(a,\nabla)v_{a}-\nabla
p_{a}+f(t,x)\right]  dt+\sum_{r=1}^{q}\gamma_{r}(t,x)dw_{r}(t), \label{os1}\\
\ \ 0\leq t\leq T,\ x\in\mathbf{R}^{n}, \nonumber \\
\operatorname{div}\ v_{a}=0, \label{os2}%
\end{gather}
with spatial periodic conditions%
\begin{align}
v_{a}(t,x+Le_{i})  &  =v_{a}(t,x),\ p_{a}(t,x+Le_{i})=p_{a}(t,x),\ \label{0s3}%
\\
0  &  \leq t\leq T,\ i=1,\ldots,n,\nonumber
\end{align}
and the initial condition%
\begin{equation}
v_{a}(0,x)=\varphi(x), \label{os4}%
\end{equation}
where $a=a(t,x)$ is an $n$-dimensional vector $a=(a^{1},\ldots,a^{n}%
)^{\intercal}$ with $a^{i}$ being $Q$-periodic deterministic functions which
have continuous derivatives with respect to $x$ up to some order; and the rest
of the notation is the same as in (\ref{NS1})-(\ref{NS3}).

We re-write the problem (\ref{os1})-(\ref{os4}) with positive direction of
time into the problem with negative direction of time which is more convenient
for making use of probabilistic representations. To this end, introduce the
new time variable $s=T-t$ and the functions $u_{a}(s,x):=v_{a}(T-s,x),$
$\tilde{a}(s,x):=a(T-s,x),$ $\tilde{f}(s,x):=f(T-s,x),$ $\tilde{\gamma}%
_{r}(s,x):=\gamma_{r}(T-s,x),$ and $\tilde{p}_{a}(s,x):=p_{a}(T-s,x).$

Further, we recall the definition of a backward Ito integral \cite{R}.
Introduce the \textquotedblleft backward\textquotedblright\ Wiener processes
\begin{equation}
\tilde{w}_{r}(t):=w_{r}(T)-w_{r}(T-t),\ \ r=1,\ldots,q,\ \ 0\leq t\leq T,
\label{bs4}%
\end{equation}
and a decreasing family of $\sigma$-subalgebras $\mathcal{F}_{t,T}^{w},$
$0\leq t\leq T,$ induced by the increments $w_{r}(T)-w_{r}(t^{\prime}),$
$r=1,\ldots,q,$ $t^{\prime}\geq t$. The increasing family of $\sigma$-subalgebras
$\mathcal{F}_{t}^{\tilde{w}}$ induced by $\tilde{w}_{r}(s^{\prime}),$
$s^{\prime}\leq t,$ coincides with $\mathcal{F}_{T-t,T}^{w},$ while
$\mathcal{F}_{t,T}^{\tilde{w}}$ is induced by the increments $\tilde{w}%
_{r}(T)-\tilde{w}_{r}(t^{\prime}),$ $r=1,\ldots,q,$ $t^{\prime}\geq t$, and
coincides with $\mathcal{F}_{T-t}^{w}.$ The backward Ito integral with respect
to $\tilde{w}_{r}(s)$ is defined as the Ito integral with respect to
$w_{r}(s)$:
\begin{equation}
\int_{t}^{t^{\prime}}\psi(t^{\prime\prime})\ast d\tilde{w}_{r}(t^{\prime
\prime}):=\int_{T-t^{\prime}}^{T-t}\psi(T-t^{\prime\prime})dw_{r}%
(t^{\prime\prime}),\ \ 0\leq t\leq t^{\prime}\leq T, \label{bs5}%
\end{equation}
where $\psi(T-t),$ $t\leq T,$ is an $\mathcal{F}_{t}^{w}$-adapted
square-integrable function and $\psi(t)$ is $\mathcal{F}_{t}^{\tilde{w}}%
$-adapted. Note that $w_{r}(t)=\tilde{w}_{r}(T)-\tilde{w}_{r}(T-t),$%
\ $r=1,\ldots,q,$\ $0\leq t\leq T.$

The backward stochastic Oseen-Stokes equations can be written as
\begin{gather}
-du_{a}(s)=\left[  \frac{\sigma^{2}}{2}\Delta u_{a}-(\tilde{a},\nabla
)u_{a}-\nabla\tilde{p}_{a}+\tilde{f}(s,x)\right]  ds+\sum_{r=1}^{q}%
\tilde{\gamma}_{r}(s,x)\ast d\tilde{w}_{r}(s),\ \ \label{os11}\\
0\leq s\leq T,\ x\in\mathbf{R}^{n},\nonumber\\
\operatorname{div}\ u_{a}=0, \label{os12}%
\end{gather}
with spatial periodic conditions%
\begin{align}
u_{a}(s,x+Le_{i})  &  =u_{a}(s,x),\ \tilde{p}_{a}(s,x+Le_{i})=\tilde{p}%
_{a}(s,x),\ \label{os13}\\
0  &  \leq s\leq T,\ i=1,\ldots,n,\nonumber
\end{align}
and the terminal condition%
\begin{equation}
u_{a}(T,x)=\varphi(x). \label{os14}%
\end{equation}
We note that (\ref{bs5}) implies
\[
\int_{s}^{T}\tilde{\gamma}_{r}(s^{\prime},x)\ast d\tilde{w}_{r}(s^{\prime
})=\int_{0}^{T-s}\gamma_{r}(s^{\prime},x)dw_{r}(s^{\prime}).
\]
The processes $u_{a}(s,x),$ $\tilde{p}_{a}(s,x)$ are $\mathcal{F}%
_{s,T}^{\tilde{w}}$-adapted (and $\mathcal{F}_{T-s}^{w}$-adapted), they depend
on $\tilde{w}_{r}(T)-\tilde{w}_{r}(s^{\prime})=w_{r}(T-s^{\prime}),$ $s\leq
s^{\prime}\leq T.$

Let $u_{a}(s,x),$\ $\tilde{p}_{a}(s,x)$ be a solution of the problem
(\ref{os11})-(\ref{os14}). For the function $u_{a}(s,x)$, one can use the
following probabilistic representation of solutions to the Cauchy problem for
linear SPDE of parabolic type (the conditional Feynman-Kac formula or the
averaging over characteristics formula, see, e.g., \cite{R} and \cite{spde}):%
\begin{equation}
u_{a}(s,x)=E^{\tilde{w}}\left[  \varphi(X_{s,x}(T))Y_{s,x,1}(T)+Z_{s,x,1,0}%
(T)\right]  ,\ 0\leq s\leq T, \label{FBD5}%
\end{equation}
where $X_{s,x}(s^{\prime}),\ Y_{s,x,y}(s^{\prime}),\ Z_{s,x,y,z}(s^{\prime
}),\ s^{\prime}\geq s,$ solves the system of Ito stochastic differential
equations:%
\begin{gather}
dX=(-\tilde{a}(s^{\prime},X)-\sigma\mu(s^{\prime},X))ds^{\prime}+\sigma
dW(s^{\prime}),\ X(s)=x,\label{BDF0}\\
dY=\mu^{\intercal}(s^{\prime},X)YdW(s^{\prime}),\ Y(s)=y,\label{BDF1}\\
dZ=(-\nabla\tilde{p}_{a}(s^{\prime},X)+\tilde{f}(s^{\prime},X))Yds^{\prime
}+F(s^{\prime},X)YdW(s^{\prime})\label{BDF2}\\
+\sum_{r=1}^{q}\tilde{\gamma}_{r}(s^{\prime},X)Yd\tilde{w}_{r}(s^{\prime
}),\ Z(s)=z.\nonumber
\end{gather}
In (\ref{FBD5})-(\ref{BDF2}), $W(s)$ is a standard $n$-dimensional Wiener
process independent of $\tilde{w}_{r}(s)$ on the probability space
$(\Omega,\mathcal{F},P)$; $Y$ is a scalar, and $Z$ is an $n$-dimensional
column-vector;$\ \mu(s,x)$ is an arbitrary $n$-dimensional spatial periodic
vector function and $F(s,x)$ is an arbitrary $n\times n$-dimensional spatial
periodic matrix function, which are sufficiently smooth in $s,x$; the
expectation $E^{\tilde{w}}$ in (\ref{FBD5}) is taken over the realizations of
$W(s),$ $t\leq s\leq T,$ for a fixed $\tilde{w}_{r}(s^{\prime}),$
$r=1,\ldots,q,$ $s\leq s^{\prime}\leq T,$ in other words, $E^{\tilde{w}%
}\left(  \cdot\right)  $ means the conditional expectation:
\[
E\left(  \cdot|\tilde{w}_{r}(s^{\prime})-\tilde{w}_{r}(s),\text{ }%
r=1,\ldots,q,\text{ }s\leq s^{\prime}\leq T\right)  .
\]

The probabilistic representation like (\ref{FBD5})-(\ref{BDF2}) for the Cauchy
problem (\ref{os11}), (\ref{os14}) is obtained (see, e.g., \cite{R}) for
linear SPDEs with deterministic coefficients. However here $\tilde{p}%
_{a}(s,x)$ is a part of solution of problem (\ref{os11})-(\ref{os14}) and it
is random (more precisely it is $\mathcal{F}_{s,T}^{\tilde{w}}$-adapted). In
this case the representation (\ref{FBD5})-(\ref{BDF2}) can be rigorously
justified in the following way. The solution $u_{a}$ of (\ref{os11}),
(\ref{os14}) can be represented in the form of the sum%
\[
u_{a}=u_{a}^{(0)}+u_{a}^{(1)},
\]
where $u_{a}^{(0)}$ satisfies the Cauchy problem for the backward
deterministic linear parabolic PDE with random parameters:%
\begin{align}
-\frac{\partial u_{a}^{(0)}}{\partial s}  &  =\frac{\sigma^{2}}{2}\Delta
u_{a}^{(0)}-(\tilde{a},\nabla)u_{a}^{(0)}-\nabla\tilde{p}_{a},\label{Fr1}\\
u_{a}^{(0)}(T,x)  &  =0,\nonumber
\end{align}
and $u_{a}^{(1)}$ satisfies the Cauchy problem for the backward stochastic
linear parabolic PDE with deterministic parameters:%
\begin{align}
-du_{a}^{(1)}(s)  &  =\left[  \frac{\sigma^{2}}{2}\Delta u_{a}^{(1)}%
-(\tilde{a},\nabla)u_{a}^{(1)}+\tilde{f}(s,x)\right]  ds+\sum_{r=1}^{q}%
\tilde{\gamma}_{r}(s,x)\ast d\tilde{w}_{r}(s),\ \label{Fr2}\\
u_{a}^{(1)}(T,x)  &  =\varphi(x).\ \nonumber
\end{align}

Clearly,%
\[
u_{a}^{(0)}(s,x)=E^{\tilde{w}}\left[  Z_{s,x,1,0}^{(0)}(T)\right]
=-E^{\tilde{w}}\int_{s}^{T}\nabla\tilde{p}_{a}(s^{\prime},X_{s,x}(s^{\prime
}))\ Y_{s,x,1}(s^{\prime})ds^{\prime}.
\]
The Feynman-Kac formula for $u_{a}^{(1)}$ coincides with (\ref{FBD5}%
)-(\ref{BDF2}) under $\nabla\tilde{p}_{a}(s,x)$ $=$ $0$.

Let $\mathcal{F}_{s,t}^{W}$ be a $\sigma$-algebra\ induced by $W_{r}%
(s^{\prime})-W_{r}(s),$ $r=1,\ldots,n,\ s\leq s^{\prime}\leq t.$ We note that
$\nabla\tilde{p}_{a}(s^{\prime},X_{s,x}(s^{\prime}))$ in (\ref{BDF2}) is
$\mathcal{F}_{s,s^{\prime}}^{W}\vee\mathcal{F}_{s^{\prime},T}^{\tilde{w}}%
$-adapted, where the family of $\sigma$-algebras $\mathcal{F}_{s,s^{\prime}%
}^{W}\vee\mathcal{F}_{s^{\prime},T}^{\tilde{w}}$ is neither increasing nor
decreasing in $s^{\prime}$. Consequently, $Z_{s,x,y,z}(s^{\prime})$ is
measurable with respect to $\mathcal{F}_{s,s^{\prime}}^{W}\vee\mathcal{F}%
_{s^{\prime},T}^{\tilde{w}}$ for every $s^{\prime}\in\lbrack s,T].$ Since
$\tilde{\gamma}_{r}(s^{\prime},X_{s,x}(s^{\prime}))Y(s^{\prime})$ are
independent of $\tilde{w}_{r},$ the Ito integral in (\ref{BDF2}) is well defined.

\begin{remark}
\label{Rem_ant}We remark that within the non-anticipating stochastic calculus
the probabilistic representation $(\ref{FBD5})$-$(\ref{BDF2})$ for the linear
problem $(\ref{os11})$-$(\ref{os14})$ cannot be carried over to the backward
SNSE problem by changing the coefficient $\tilde{a}(s,x)$ to $u(s,x)$ since
then the integrand $\tilde{\gamma}_{r}(s^{\prime},X_{s,x}(s^{\prime
}))Y(s^{\prime})$ would be $\mathcal{F}_{s,s^{\prime}}^{W}\vee\mathcal{F}%
_{s^{\prime},T}^{\tilde{w}}$-measurable. Nevertheless, the representation
$(\ref{FBD5})$-$(\ref{BDF2})$ allows us to derive layer methods for the
stochastic Oseen-Stokes equations $(\ref{os11})$-$(\ref{os14})$, and then,
using them as a guidance, one can obtain layer methods for the SNSE
$(\ref{NS1})$-$(\ref{NS3})$ as well $($see Sections~3.1 and~3.2$)$.
\end{remark}

For deriving layer methods, we also use some direct probabilistic
representations for solutions of the SNSE. In Sections~\ref{prepDirect}
and~\ref{prepDoub} we give two such representations. The first one follows
from a specific probabilistic representation for a linear SPDE which differs
from $(\ref{FBD5})$-$(\ref{BDF2})$ and the second one uses backward doubly
stochastic differential equations \cite{PP}.

\subsection{A direct probabilistic representation for solutions of
SNSE\label{prepDirect}}

As in the case of the stochastic Oseen-Stokes equations, we re-write the SNSE
problem (\ref{NS1})-(\ref{NS3}) with positive direction of time into the
problem with negative direction of time. Again introduce the new time variable
$s=T-t$ and the functions $u(s,x):=v(T-s,x)$, $\tilde{f}(s,x):=f(T-s,x),$
$\tilde{\gamma}_{r}(s,x):=\gamma_{r}(T-s,x),$ and $\tilde{p}(s,x):=p(T-s,x).$
The corresponding backward SNSE take the form:
\begin{gather}
-du=(\frac{\sigma^{2}}{2}\Delta u-(u,\nabla)u-\nabla\tilde{p}+\tilde
{f})ds+\sum_{r=1}^{q}\tilde{\gamma}_{r}(s,x)\ast d\tilde{w}_{r}%
(s),\ u(T,x)=\varphi(x),\label{D1}\\
\operatorname{div}u=0\ , \label{D01}%
\end{gather}
with spatial periodic conditions for $u$ and $\tilde{p}$.

Introduce $F(s,x,u,\nabla u):=-(u,\nabla)u-\nabla\tilde{p}+\tilde{f}$ and
write (\ref{D1}) as
\begin{equation}
-du=\left(  \frac{\sigma^{2}}{2}\Delta u+F(s,x,u,\nabla u)\right)
ds+\sum_{r=1}^{q}\tilde{\gamma}_{r}(s,x)\ast d\tilde{w}_{r}%
(s),\ u(T,x)=\varphi(x). \label{D03}%
\end{equation}
Let us assume that the solution $u(s,x)=u(s,x,\omega)$ to (\ref{D1}%
)-(\ref{D01}) is known. We substitute it in the $F(s,x,u,\nabla u)$ which
becomes a function $\tilde{F}(s,x,\omega)$ depending on $\omega$ as a
parameter. Hence (\ref{D03}) can be considered as a linear parabolic SPDE. For
solutions of this linear SPDE, we can write the following probabilistic
representation analogously to (\ref{FBD5})-(\ref{BDF2}) (we take $Y\equiv1$):
\begin{gather}
u(s,x)=E^{\tilde{w}}\varphi(X_{s,x}(T))\label{NS11n}\\
-E^{\tilde{w}}\left[  \int_{s}^{T}\{\nabla\tilde{p}(s^{\prime},X_{s,x}%
(s^{\prime}))-\tilde{f}(s^{\prime},X_{s,x}(s^{\prime})) \right. \nonumber \\
\Bigg. +(u(s^{\prime}%
,X_{s,x}(s^{\prime})),\nabla)u(s^{\prime},X_{s,x}(s^{\prime}))\}ds^{\prime
}\Bigg] \nonumber\\
+\sum_{r=1}^{q}E^{\tilde{w}}\left[  \int_{s}^{T}\tilde{\gamma}_{r}(s^{\prime
},X_{s,x}(s^{\prime}))d\tilde{w}_{r}(s^{\prime})\right]  ,\nonumber
\end{gather}
where $X_{s,x}(s^{\prime}),$ $s^{\prime}\geq s,$ solves the system of
stochastic differential equations
\begin{equation}
dX=\sigma dW(s^{\prime}),\ X(s)=x, \label{NS12n}%
\end{equation}
$W$ is a standard $n$-dimensional Wiener process independent of $\tilde{w}%
_{r}$ on the probability space $(\Omega,\mathcal{F},P).$

\subsection{A probabilistic representation for solution of SNSE using backward
doubly stochastic differential equations\label{prepDoub}}

In connection with the backward SNSE (\ref{D1})-(\ref{D01}), we introduce the
system of backward doubly stochastic differential equations \cite{PP}:
\begin{align}
dX  &  =\sigma dW(s^{\prime}),\ \ X(s)=x,\label{D2} \displaybreak[0]\\
dU  &  =(\nabla\tilde{p}(s^{\prime},X)-\tilde{f}(s^{\prime},X)+\frac{1}%
{\sigma}\mathbb{Z}U)ds^{\prime}+\mathbb{Z}dW(s^{\prime})-\sum_{r=1}^{q}%
\tilde{\gamma}_{r}(s^{\prime},X)\ast d\tilde{w}_{r}(s^{\prime}),\label{D3} \displaybreak[0]\\
&U(T)    =\varphi(X_{s,x}(T)). \label{D4}%
\end{align}
In (\ref{D2})-(\ref{D4}) $X,$\ $U,$\ $W$ are column vectors of dimension $n$
and $\mathbb{Z}$ is a matrix of dimension $n\times n,\ W(s)$ and $\tilde
{w}(s),\ 0\leq s\leq T,$ are mutually independent standard Wiener processes on
the probability space $(\Omega,\mathcal{F},P)$. We recall that the triple
$\{X_{s,x}(s^{\prime}),U_{s,x}(s^{\prime}),\mathbb{Z}_{s,x}(s^{\prime}),s\leq
s^{\prime}\leq T\}$ is a solution of (\ref{D2})-(\ref{D4}) if $X_{s,x}%
(s^{\prime})$ satisfies (\ref{D2}), $(U_{s,x}(s^{\prime}),\mathbb{Z}%
_{s,x}(s^{\prime}))$ for each $s^{\prime}$ is $\mathcal{F}_{s,s^{\prime}}%
^{W}\vee\mathcal{F}_{s^{\prime},T}^{\tilde{w}}$-measurable, and%
\begin{gather}
U_{s,x}(s^{\prime})=\varphi(X_{s,x}(T))-\int_{s^{\prime}}^{T}(\nabla\tilde
{p}(s^{\prime\prime},X_{s,x}(s^{\prime\prime}))-\tilde{f}(s^{\prime\prime
},X_{s,x}(s^{\prime\prime})) \nonumber \\
+\frac{1}{\sigma}\mathbb{Z}_{s,x}(s^{\prime\prime
})U_{s,x}(s^{\prime\prime}))ds^{\prime\prime}\label{D45}\\
-\int_{s^{\prime}}^{T}\mathbb{Z}_{s,x}(s^{\prime\prime})dW(s^{\prime\prime
})+\int_{s^{\prime}}^{T}\sum_{r=1}^{q}\tilde{\gamma}_{r}(s^{\prime\prime
},X_{s,x}(s^{\prime\prime}))\ast d\tilde{w}_{r}(s^{\prime\prime}),\ s\leq
s^{\prime}\leq T.\nonumber
\end{gather}

Let $u(s,x)$ be a solution of the problem (\ref{D1}), i.e.,
\begin{align}
u(s,x)  &  =\varphi(x)+\int_{s}^{T}(\frac{\sigma^{2}}{2}\Delta u(s^{\prime
},x)-(u,\nabla)u(s^{\prime},x)-\nabla\tilde{p}(s^{\prime},x)+\tilde
{f}(s^{\prime},x))ds^{\prime}\label{D5}\\
&  +\sum_{r=1}^{q}\int_{s}^{T}\tilde{\gamma}_{r}(s^{\prime},x)\ast d\tilde
{w}_{r}(s^{\prime}).\nonumber
\end{align}
It is known (see \cite{PP}) that then%
\begin{align}
X(s^{\prime})  &  =X_{s,x}(s^{\prime}),\ U(s^{\prime})=U_{s,x}(s^{\prime
})=u(s^{\prime},X_{s,x}(s^{\prime})),\label{D55}\\
\mathbb{Z}(s^{\prime})  &  =\mathbb{Z}_{s,x}(s^{\prime})=\{\mathbb{Z}%
^{k,j}(s^{\prime})\}=\sigma\cdot\left\{  \frac{\partial u^{k}}{\partial x^{j}%
}(s^{\prime},X_{s,x}(s^{\prime}))\right\}  ,\ k,j=1,\ldots,n,\nonumber
\end{align}
is a solution of (\ref{D2})-(\ref{D4}).

Conversely, if $X_{s,x}(s^{\prime}),$\ $U_{s,x}(s^{\prime}),$\ $\mathbb{Z}%
_{s,x}(s^{\prime})$ is a solution of the system of backward doubly stochastic
differential equations (\ref{D2})-(\ref{D4}) then it can be verified that
\begin{equation}
u(s,x)=U_{s,x}(s) \label{D9}%
\end{equation}
is the solution of (\ref{D1}) (see \cite{PP}). The condition (\ref{D01}) is
satisfied by choosing an appropriate pressure $\tilde{p}$.

We note that $u(s,x)$ is $\mathcal{F}_{s,T}^{\tilde{w}}$-measurable and then
using (\ref{D45}) we get
\begin{align}
u(s,x)  &  =U_{s,x}(s)=E[U_{s,x}(s)|\mathcal{F}_{s,T}^{\tilde{w}}%
]=E^{\tilde{w}}U_{s,x}(s)\label{D10} \displaybreak[0]\\
&  =E^{\tilde{w}}\varphi(X_{s,x}(T))\nonumber \displaybreak[0]\\
&  -E^{\tilde{w}}\int_{s}^{T}(\nabla\tilde{p}(s^{\prime},X_{s,x}(s^{\prime
}))-\tilde{f}(s^{\prime},X_{s,x}(s^{\prime}))+\frac{1}{\sigma}\mathbb{Z}%
_{s,x}(s^{\prime})U_{s,x}(s^{\prime}))ds^{\prime}\nonumber \displaybreak[0]\\
&  +\sum_{r=1}^{q}E^{\tilde{w}}\int_{s}^{T}\tilde{\gamma}_{r}(s^{\prime
},X_{s,x}(s^{\prime}))\ast d\tilde{w}_{r}(s^{\prime}).\nonumber
\end{align}
Due to smoothness of $\tilde{\gamma}_{r}(s,x)$ in $s$ and independence of $X$
and $\tilde{w},$ the equality%
\[
\int_{s}^{T}\tilde{\gamma}_{r}(s^{\prime},X_{s,x}(s^{\prime}))\ast d\tilde
{w}_{r}(s^{\prime})=\int_{s}^{T}\tilde{\gamma}_{r}(s^{\prime},X_{s,x}%
(s^{\prime}))d\tilde{w}_{r}(s^{\prime})
\]
holds. Hence the right-hand side of (\ref{D10}) coincides with the right-hand
side of the probabilistic representation (\ref{NS11n}).

\section{Layer methods\label{secLayer}}

In this section we construct three layer methods based on the probabilistic
representations from Sections~\ref{prepOseen} and~\ref{prepDirect}. In the
case of deterministic NSE (i.e., when $\gamma_{r}=0$ in the SNSE
(\ref{NS1})-(\ref{NS3})) these methods coincide with the ones presented in
\cite{NS5}.

On the basis of the probabilistic representation (\ref{FBD5})-(\ref{BDF2}) we,
first, construct layer methods for the stochastic Oseen-Stokes equations and,
second, using the obtained methods as a guidance, we construct the
corresponding methods for the SNSE (this way of deriving numerical methods for
nonlinear SPDEs was proposed in \cite{spde}). This is done in
Sections~\ref{secSL} and~\ref{secSimL}. We underline that derivation of these
methods does not rely on direct probabilistic representations for the SNSE
themselves that would require the anticipating stochastic calculus (see
Remark~\ref{Rem_ant}) which is not developed satisfactorily from the numerical
point of view. That is why we prefer to use the mimicry approach here.

In Section~\ref{secDiL}\ we derive a layer method based on the direct
probabilistic representation for the SNSE from Section~\ref{prepDirect}.

In Sections~\ref{secSL}, \ref{secSimL} and~\ref{secDiL} we deal with
approximation of velocity $v(t,x)$ (i.e., a part of the solution $v(t,x),$
$p(t,x)$ to the SNSE) only. Since we consider here the spatial-periodic
problem (\ref{NS1})-(\ref{NS3}), we can separate approximation of velocity
$v(t,x)$ and pressure $p(t,x)$ in a constructive way. Approximation of
pressure is considered in Section~\ref{secPres}.

Let us introduce a uniform partition of the time interval $[0,T]:$
$0=t_{0}<t_{1}<\cdots<t_{N}=T$ and the time step $h=T/N$ (we restrict
ourselves to the uniform partition for simplicity only).

\subsection{A layer method based on the standard probabilistic
representation\label{secSL}}

Each choice of $\mu(s,x)$ and $F(s,x)$ in (\ref{FBD5})-(\ref{BDF2}) gives us a
particular probabilistic representation for the solution of the stochastic
Oseen-Stokes equations (\ref{os11})-(\ref{os14}) which can be used for
deriving the corresponding layer method. In this and the next section we
derive layer methods based on two of such probabilistic representations which
can be, in a sense, viewed as limiting cases of (\ref{FBD5})-(\ref{BDF2}). If
we put $\mu(s,x)=0$ and $F(s,x)=0$ in (\ref{FBD5})-(\ref{BDF2}), we obtain the
standard probabilistic representation for the solution to the backward linear
SPDE (\ref{os11})-(\ref{os14}) \cite{R}. This case is considered in this
section. The case of $F(s,x)=0$ and $\mu(s,x)$ turning the equation
(\ref{BDF0}) for $X(s)$ into pure diffusion is treated in the next section.

Analogously to (\ref{FBD5})-(\ref{BDF2}) with $\mu(s,x)=0$ and $F(s,x)=0,$ we
get the following local probabilistic representation of the solution to
(\ref{os11})-(\ref{os14}):
\begin{align}
u_{a}(t_{k},x)  &  =E^{\tilde{w}}\left[  u_{a}(t_{k+1},X_{t_{k},x}%
(t_{k+1}))-\int_{t_{k}}^{t_{k+1}}\nabla\tilde{p}_{a}(s,X_{t_{k},x}%
(s))ds\right. \label{spr1}\\
&  \left.  +\int_{t_{k}}^{t_{k+1}}\tilde{f}(s,X_{t_{k},x}(s))ds+\sum_{r=1}%
^{q}\int_{t_{k}}^{t_{k+1}}\tilde{\gamma}_{r}(s,X_{t_{k},x}(s))d\tilde{w}%
_{r}(s)\right]  ,\nonumber
\end{align}
where
\begin{equation}
dX=-\tilde{a}(s,X)ds+\sigma dW(s),\ X(t_{k})=x. \label{spr2}%
\end{equation}

A slightly modified explicit Euler scheme with the simplest noise simulation
applied to (\ref{spr2}) gives%
\begin{equation}
X_{t_{k},x}(t_{k+1})\simeq\bar{X}_{t_{k},x}(t_{k+1})=x-\tilde{a}%
(t_{k+1},x)h+\sigma\sqrt{h}\xi, \label{NS15}%
\end{equation}
where $\xi=(\xi^{1},\ldots,\xi^{n})^{\top}$ and$\ \xi^{1},\ldots,\xi^{n}$ are
i.i.d. random variables with the law $P(\xi^{i}=\pm1)=1/2.$ We substitute
$\bar{X}_{t_{k},x}(t_{k+1})$ from (\ref{NS15}) in (\ref{spr1}) instead of
$X_{t_{k},x}(t_{k+1})$, evaluate the expectation exactly, and thus obtain
(recall that $\operatorname{div}\tilde{\gamma}_{r}=0$ and $\nabla\tilde{p}%
_{a}(s,x)\in(\mathbf{V}_{p}^{m})^{\bot}):$
\begin{gather}
u_{a}(t_{k},x)=\breve{u}_{a}(t_{k+1},x)-\nabla\tilde{p}_{a}(t_{k+1}%
,x)h+\tilde{f}(t_{k+1},x)h\label{NS16}\\
+\sum_{r=1}^{q}\tilde{\gamma}_{r}(t_{k+1},x)\left(  \tilde{w}_{r}%
(t_{k+1})-\tilde{w}_{r}(t_{k})\right)  +\rho\nonumber\\
=P\breve{u}_{a}(t_{k+1},x)+P\tilde{f}(t_{k+1},x)h+P^{\bot}\breve{u}%
_{a}(t_{k+1},x)+P^{\bot}\tilde{f}(t_{k+1},x)h\nonumber\\
-\nabla\tilde{p}_{a}%
(t_{k+1},x)h
+\sum_{r=1}^{q}\tilde{\gamma}_{r}(t_{k+1},x)\Delta_{k}\tilde{w}_{r}%
+\rho,\nonumber
\end{gather}
where $\Delta_{k}\tilde{w}_{r}=\tilde{w}_{r}(t_{k+1})-\tilde{w}_{r}(t_{k}),$
$r=1,\ldots,q;$ $\rho=\rho(t_{k},x)$ is a remainder, and
\begin{equation}
\breve{u}_{a}(t_{k+1},x)=E^{\tilde{w}}u_{a}(t_{k+1},\bar{X}_{k+1})=2^{-n}%
\sum_{j=1}^{2^{n}}u_{a}(t_{k+1},x-\tilde{a}(t_{k+1},x)h+\sigma\sqrt{h}\xi_{j})
\label{NS17}%
\end{equation}
with $\xi_{1}=(1,1,\ldots,1)^{\top},\ \ldots,\ \xi_{2^{n}}=(-1,-1,\ldots
,-1)^{\top}.$ Taking into account that $u_{a}(t_{k},x)$ in (\ref{NS16}) is
divergence free, we get%
\begin{equation}
u_{a}(t_{k},x)=P\breve{u}_{a}(t_{k+1},x)+P\tilde{f}(t_{k+1},x)h+\sum_{r=1}%
^{q}\tilde{\gamma}_{r}(t_{k+1},x)\Delta_{k}\tilde{w}_{r}+P\rho. \label{NS165}%
\end{equation}
Neglecting the remainder, we get the one-step approximation for $u_{a}%
(t_{k},x)$:
\begin{equation}
\hat{u}_{a}(t_{k},x)=P\breve{u}_{a}(t_{k+1},x)+P\tilde{f}(t_{k+1}%
,x)h+\sum_{r=1}^{q}\tilde{\gamma}_{r}(t_{k+1},x)\Delta_{k}\tilde{w}_{r}.
\label{os15n}%
\end{equation}

Re-writing $\hat{u}_{a}(t_{k},x)$ of (\ref{os15n}) in the positive direction
of time, we obtain the one-step approximation for the velocity $v_{a}%
(t_{k},x)$ of the forward-time stochastic Oseen-Stokes equations
(\ref{os1})-(\ref{os4}):
\begin{equation}
\hat{v}_{a}(t_{k+1},x)=P\breve{v}_{a}(t_{k},x)+Pf(t_{k},x)h+\sum_{r=1}%
^{q}\gamma_{r}(t_{k},x)\Delta_{k}w_{r}, \label{os15}%
\end{equation}
where $\Delta_{k}w_{r}=w_{r}(t_{k+1})-w_{r}(t_{k}),$ $r=1,\ldots,q,$ and
\begin{equation}
\breve{v}_{a}(t_{k},x)=2^{-n}\sum_{j=1}^{2^{n}}v_{a}(t_{k},x-a(t_{k}%
,x)h+\sigma\sqrt{h}\xi_{j}). \label{os17}%
\end{equation}

Now let us turn our attention from the stochastic Oseen-Stokes equation to the
\textit{stochastic NSE} (\ref{NS1})-(\ref{NS3}).

Using the one-step approximation (\ref{os15})-(\ref{os17}) for the stochastic
Oseen-Stokes equations (\ref{os1})-(\ref{os4}) as a guidance, we construct the
one-step approximation for the SNSE (\ref{NS1})-(\ref{NS3}) by substituting
$a(t_{k},x)$ with $v(t_{k},x):$
\begin{equation}
\hat{v}(t_{k+1},x)=P\breve{v}(t_{k},x)+Pf(t_{k},x)h+\sum_{r=1}^{q}\gamma
_{r}(t_{k},x)\Delta_{k}w_{r}, \label{NSA1}%
\end{equation}
where
\begin{equation}
\breve{v}(t_{k},x)=2^{-n}\sum_{j=1}^{2^{n}}v(t_{k},x-v(t_{k},x)h+\sigma
\sqrt{h}\xi_{j}). \label{NSA3}%
\end{equation}
It is easy to see that under Assumptions~2.1 $\operatorname{div}\hat
{v}(t_{k+1},x)=0.$

The corresponding layer method for the SNSE (\ref{NS1})-(\ref{NS3}) has the
form%
\begin{gather}
\bar{v}(0,x)=\varphi(x),\ \bar{v}(t_{k+1},x)=P\breve{v}(t_{k},x)+Pf(t_{k}%
,x)h+\sum_{r=1}^{q}\gamma_{r}(t_{k},x)\Delta_{k}w_{r},\label{NS18}\\
k=0,\ldots,N-1,\nonumber
\end{gather}
where%
\begin{equation}
\breve{v}(t_{k},x)=2^{-n}\sum_{j=1}^{2^{n}}\bar{v}(t_{k},x-\bar{v}%
(t_{k},x)h+\sigma\sqrt{h}\xi_{j}). \label{NS19}%
\end{equation}
We note that we use the same notation $\breve{v}(t_{k},x)$ for the functions
appearing in the one-step approximation (\ref{NSA3}) and in the layer method
(\ref{NS19}) but this does not cause any confusion.

Knowing the expansions
\begin{align}
\breve{v}(t_{k},x)  &  =\sum_{\mathbf{n}\in\mathbf{Z}^{n}}\breve
{v}_{\mathbf{n}}(t_{k})e^{i(2\pi/L)(\mathbf{n},x)},\ \ \ f(t_{k}%
,x)=\sum_{\mathbf{n}\in\mathbf{Z}^{n}}f_{\mathbf{n}}(t_{k})e^{i(2\pi
/L)(\mathbf{n},x)},\label{NS20}\\
\gamma_{r}(t_{k},x)  &  =\sum_{\mathbf{n}\in\mathbf{Z}^{n}}\gamma
_{r,\mathbf{n}}(t_{k})e^{i(2\pi/L)(\mathbf{n},x)},\nonumber
\end{align}
it is not difficult to find $\bar{v}(t_{k+1},x)$. Indeed, using (\ref{N00})
and (\ref{N01}), we obtain from (\ref{NS18})-(\ref{NS19}):
\begin{gather}
\bar{v}(t_{k+1},x)=\sum_{\mathbf{n}\in\mathbf{Z}^{n}}\bar{v}_{\mathbf{n}%
}(t_{k+1})e^{i(2\pi/L)(\mathbf{n},x)},\ \ \label{NS21}\\
\bar{v}_{\mathbf{n}}(t_{k+1})=\breve{v}_{\mathbf{n}}(t_{k})+f_{\mathbf{n}%
}(t_{k})h-\frac{\breve{v}_{\mathbf{n}}^{\top}(t_{k})\mathbf{n}}{|\mathbf{n}%
|^{2}}\mathbf{n}-h\frac{f_{\mathbf{n}}^{\top}(t_{k})\mathbf{n}}{|\mathbf{n}%
|^{2}}\mathbf{n}+\sum_{r=1}^{q}\gamma_{r,\mathbf{n}}(t_{k})\ \Delta_{k}%
w_{r}\mathbf{.}\nonumber
\end{gather}
We note that turning the layer method (\ref{NS18})-(\ref{NS19}) into a
numerical algorithm requires to complement it with an interpolation in order
to compute the terms $\bar{v}(t_{k},x-\bar{v}(t_{k},x)h+\sigma\sqrt{h}\xi
_{j})$ in (\ref{NS19}) used for finding $\breve{v}_{\mathbf{n}}(t_{k})$ from
(\ref{NS20}), see the corresponding discussion in the case of deterministic
NSE in \cite{NS5}.

\subsection{Layer methods based on the probabilistic representation with
simplest characteristics\label{secSimL}}

If we put $\mu(s,x)=-\tilde{a}(s,x)/\sigma$ and $F(s,x)=0$ in (\ref{FBD5}%
)-(\ref{BDF2}), we can obtain the following local probabilistic representation
for the solution to the backward stochastic Oseen-Stokes equation
(\ref{os11})-(\ref{os14}):%
\begin{gather}
u_{a}(t_{k},x)=E^{\tilde{w}}[u_{a}(t_{k+1},X_{t_{k},x}(t_{k+1}))Y_{t_{k}%
,x,1}(t_{k+1})]\label{NS14}\\
+E^{\tilde{w}}\left[  -\int_{t_{k}}^{t_{k+1}}\nabla\tilde{p}_{a}(s,X_{t_{k}%
,x}(s))Y_{t_{k},x,1}(s)ds+\int_{t_{k}}^{t_{k+1}}\tilde{f}(s,X_{t_{k}%
,x}(s))Y_{t_{k},x,1}(s)ds\right. \nonumber\\
\left.  +\sum_{r=1}^{q}\int_{t_{k}}^{t_{k+1}}\tilde{\gamma}_{r}(s,X_{t_{k}%
,x}(s))Y_{t_{k},x,1}(s)d\tilde{w}_{r}(s)\right]  \ ,\nonumber
\end{gather}
where $X_{t,x}(s),$\ $Y_{t,x,1}(s),$ $s\geq t,$ solve the system of stochastic
differential equations
\begin{align}
dX  &  =\sigma dW(s),\ X(t)=x,\label{NS12}\\
dY  &  =-\frac{1}{\sigma}Y\tilde{a}^{\top}(s,X)dW(s),\ Y(t)=1. \label{NS13}%
\end{align}
We apply a slightly modified explicit Euler scheme with the simplest noise
simulation to (\ref{NS12})-(\ref{NS13}):
\begin{equation}
\bar{X}_{t_{k},x}(t_{k+1})=x+\sigma\sqrt{h}\xi,\ \bar{Y}_{t_{k},x,1}%
(t_{k+1})=1-\frac{1}{\sigma}\tilde{a}^{\top}(t_{k+1},x)\sqrt{h}\xi,
\label{NS30}%
\end{equation}
where $\xi$ is the same as in (\ref{NS15}). Approximating $X_{t_{k},x}%
(t_{k+1})$ and $Y_{t_{k},x,1}(t_{k+1})$ in (\ref{NS14}) by $\bar{X}_{t_{k}%
,x}(t_{k+1})$ and $\bar{Y}_{t_{k},x,1}(t_{k+1})$ from (\ref{NS30}), we obtain%
\begin{gather}
u_{a}(t_{k},x)=E^{\tilde{w}}[u_{a}(t_{k+1},x+\sigma\sqrt{h}\xi)(1-\frac
{1}{\sigma}\tilde{a}^{\top}(t_{k+1},x)\sqrt{h}\xi)]-\nabla\tilde{p}%
_{a}(t_{k+1},x)h\label{NS31}\\
+\tilde{f}(t_{k+1},x)h+\sum_{r=1}^{q}\tilde{\gamma}_{r}(t_{k+1},x)\Delta
_{k}\tilde{w}_{r}+\rho\nonumber\\
=2^{-n}\sum_{q=1}^{2^{n}}u_{a}(t_{k+1},x+\sigma\sqrt{h}\xi_{q})-\frac{\sqrt
{h}}{\sigma}\breve{u}_{a}(t_{k+1},x)-\nabla\tilde{p}_{a}(t_{k+1},x)h \nonumber\\
+\tilde{f}(t_{k+1},x)h
+\sum_{r=1}^{q}\tilde{\gamma}_{r}(t_{k+1},x)\Delta_{k}\tilde{w}_{r}%
+\rho,\nonumber
\end{gather}
where
\begin{align}
\breve{u}_{a}(t_{k+1},x)  &  =E^{\tilde{w}}[u_{a}(t_{k+1},x+\sigma\sqrt{h}%
\xi)\xi^{\top}]\tilde{a}(t_{k+1},x)\label{NS32}\\
&  =2^{-n}\sum_{j=1}^{2^{n}}u_{a}(t_{k+1},x+\sigma\sqrt{h}\xi_{j})\xi
_{j}^{\top}\tilde{a}(t_{k+1},x)\nonumber
\end{align}
and $\rho=\rho(t_{k},x)$ is a remainder.

Using the Helmholtz-Hodge-Leray decomposition and taking into account that
\[
\operatorname{div}u_{a}(t_{k+1},x+\sigma\sqrt{h}\xi_{q})=0,\text{\ \ }%
\operatorname{div}\gamma_{r}=0,
\]
\ we get from (\ref{NS31})-(\ref{NS32}):
\begin{gather*}
u_{a}(t_{k},x)=2^{-n}\sum_{j=1}^{2^{n}}u_{a}(t_{k+1},x+\sigma\sqrt{h}\xi
_{j})-\frac{\sqrt{h}}{\sigma}P\breve{u}_{a}(t_{k+1},x)+P\tilde{f}%
(t_{k+1},x)h\\
-\frac{\sqrt{h}}{\sigma}P^{\bot}\breve{u}_{a}(t_{k+1},x)+P^{\bot}\tilde
{f}(t_{k+1},x)h-\nabla\tilde{p}_{a}(t_{k+1},x)h\\
+\sum_{r=1}^{q}\tilde{\gamma}_{r}(t_{k+1},x)\Delta_{k}\tilde{w}_{r}+\rho,
\end{gather*}
whence we obtain after applying the operator $P:$
\begin{gather}
u_{a}(t_{k},x)=2^{-n}\sum_{j=1}^{2^{n}}u_{a}(t_{k+1},x+\sigma\sqrt{h}\xi
_{j})-\frac{\sqrt{h}}{\sigma}P\breve{u}_{a}(t_{k+1},x)+P\tilde{f}%
(t_{k+1},x)h\label{NS33} \displaybreak[0]\\
+\sum_{r=1}^{q}\tilde{\gamma}_{r}(t_{k+1},x)\Delta_{k}\tilde{w}_{r}%
+P\rho.\nonumber
\end{gather}
Dropping the remainder in (\ref{NS33}) and re-writing the obtained
approximation in the one with positive direction of time, we obtain the
one-step approximation for the forward-time stochastic Oseen-Stokes equation
(\ref{os1})-(\ref{os4}):
\begin{gather}
\hat{v}_{a}(t_{k+1},x)=2^{-n}\sum_{j=1}^{2^{n}}v_{a}(t_{k},x+\sigma\sqrt{h}%
\xi_{j})-\frac{\sqrt{h}}{\sigma}P\breve{v}_{a}(t_{k},x)+Pf(t_{k}%
,x)h\label{os18}\\
+\sum_{r=1}^{q}\gamma_{r}(t_{k},x)\Delta_{k}w_{r},\nonumber
\end{gather}
where%
\begin{equation}
\breve{v}_{a}(t_{k},x)=2^{-n}\sum_{j=1}^{2^{n}}v_{a}(t_{k},x+\sigma\sqrt{h}%
\xi_{j})\xi_{j}^{\top}a(t_{k},x). \label{os20}%
\end{equation}
Using (\ref{os18})-(\ref{os20}) as a guidance, we arrive at the one-step
approximation for the SNSE (\ref{NS1})-(\ref{NS3}):
\begin{gather}
\hat{v}(t_{k+1},x)=2^{-n}\sum_{q=1}^{2^{n}}v(t_{k},x+\sigma\sqrt{h}\xi
_{q})-\frac{\sqrt{h}}{\sigma}P\breve{v}(t_{k},x) \label{NSA4} \\
+Pf(t_{k},x)h+\sum_{r=1}%
^{q}\gamma_{r}(t_{k},x)\Delta_{k}w_{r}, \nonumber
\end{gather}
where%
\begin{equation}
\breve{v}(t_{k},x)=2^{-n}\sum_{j=1}^{2^{n}}v(t_{k},x+\sigma\sqrt{h}\xi_{j}%
)\xi_{j}^{\top}v(t_{k},x). \label{NSA6}%
\end{equation}
It is easy to see that under Assumptions~2.1 $\operatorname{div}\hat
{v}(t_{k+1},x)=0.$ The corresponding layer method for the SNSE (\ref{NS1}%
)-(\ref{NS3}) has the form%
\begin{gather}
\bar{v}(0,x)=\varphi(x),\ \bar{v}(t_{k+1},x)=2^{-n}\sum_{j=1}^{2^{n}}\bar
{v}(t_{k},x+\sigma\sqrt{h}\xi_{j})-\frac{\sqrt{h}}{\sigma}P\breve{v}%
(t_{k},x)\label{NSM21} \displaybreak[0]\\
+Pf(t_{k},x)h+\sum_{r=1}^{q}\gamma_{r}(t_{k},x)\Delta_{k}w_{r},\ \ k=0,\ldots
,N-1,\nonumber
\end{gather}
where
\begin{equation}
\breve{v}(t_{k},x)=2^{-n}\sum_{j=1}^{2^{n}}\bar{v}(t_{k},x+\sigma\sqrt{h}%
\xi_{j})\xi_{j}^{\top}\bar{v}(t_{k},x). \label{NSM23}%
\end{equation}

Practical implementation of the layer method (\ref{NSM21})-(\ref{NSM23}) is
straightforward and efficient. Let us write the corresponding numerical
algorithm for simplicity in the two-dimensional ($n=2)$ case. We choose a
positive integer $M$ as a cut-off frequency and write the approximate velocity
at the time $t_{k+1}$ as the partial sum:
\begin{equation}
\bar{v}(t_{k+1},x)=\sum_{n_{1}=-M}^{M-1}\sum_{n_{2}=-M}^{M-1}\bar
{v}_{\mathbf{n}}(t_{k+1})e^{i(2\pi/L)(\mathbf{n},x)}, \label{al1}%
\end{equation}
where $\mathbf{n}=(n_{1},n_{2})^{\top}.$

We note that we use the same notation $\bar{v}(t_{k+1},x)$ for the partial sum
in (\ref{al1}) instead of writing $\bar{v}_{M}(t_{k+1},x)$ while in
(\ref{NSM21}) $\bar{v}(t_{k+1},x)$ denotes the approximate velocity containing
all frequencies but this should not lead to any confusion.

Further, we have
\begin{equation}
\frac{1}{4}\sum_{j=1}^{4}\bar{v}(t_{k},x+\sigma\sqrt{h}\xi_{j})=\sum
_{n_{1}=-M}^{M-1}\sum_{n_{2}=-M}^{M-1}\bar{v}_{\mathbf{n}}(t_{k}%
)e^{i(2\pi/L)(\mathbf{n},x)}\frac{1}{4}\sum_{j=1}^{4}e^{i(2\pi\sigma\sqrt
{h}/L)(\mathbf{n},\xi_{j})}. \label{al11}%
\end{equation}
Then
\begin{align*}
\breve{v}(t_{k},x)  &  =\frac{1}{4}\sum_{j=1}^{4}\bar{v}(t_{k},x+\sigma
\sqrt{h}\xi_{j})\xi_{j}^{\top}\bar{v}(t_{k},x)\\
&  =\sum_{n_{1}=-N}^{M-1}\sum_{n_{2}=-N}^{M-1}\bar{v}_{\mathbf{n}}%
(t_{k})e^{i(2\pi/L)(\mathbf{n},x)}\frac{1}{4}\sum_{j=1}^{4}e^{i(2\pi
\sigma\sqrt{h}/L)(\mathbf{n},\xi_{j})}\xi_{j}^{\top}\bar{v}(t_{k},x)\\
&  =\sum_{n_{1}=-M}^{M-1}\sum_{n_{2}=-M}^{M-1}V_{\mathbf{n}}(t_{k}%
)e^{i(2\pi/L)(\mathbf{n},x)}\bar{v}(t_{k},x),
\end{align*}
where
\[
V_{\mathbf{n}}(t_{k})=\bar{v}_{\mathbf{n}}(t_{k})\cdot\frac{1}{4}\sum
_{j=1}^{4}e^{i(2\pi\sigma\sqrt{h}/L)(\mathbf{n},\xi_{j})}\xi_{j}^{\top}.
\]
Note that $V_{\mathbf{n}}(t_{k})$ is a $2\times2$-matrix. Let
\begin{equation}
V(t_{k},x):=\sum_{n_{1}=-M}^{M-1}\sum_{n_{2}=-M}^{M-1}V_{\mathbf{n}}%
(t_{k})e^{i(2\pi/L)(\mathbf{n},x)} \label{ext}%
\end{equation}
then
\[
\breve{v}(t_{k},x)=V(t_{k},x)\bar{v}(t_{k},x).
\]

We obtain the algorithm:
\begin{align}
\bar{v}_{\mathbf{n}}(0)  &  =\varphi_{\mathbf{n}},\ \label{alg2} \displaybreak[0]\\
\bar{v}_{\mathbf{n}}(t_{k+1})  &  =\bar{v}_{\mathbf{n}}(t_{k})-\frac{\sqrt{h}%
}{\sigma}\left(  \breve{v}_{\mathbf{n}}(t_{k})-\frac{\breve{v}_{\mathbf{n}%
}^{\top}(t_{k})\mathbf{n}}{|\mathbf{n}|^{2}}\mathbf{n}\right)  +f_{\mathbf{n}%
}(t_{k})h-h\frac{f_{\mathbf{n}}^{\top}(t_{k})\mathbf{n}}{|\mathbf{n}|^{2}%
}\mathbf{n} \nonumber \displaybreak[0]\\
&+\sum_{r=1}^{q}\gamma_{r,\mathbf{n}}(t_{k})\ \Delta_{k}%
w_{r},\nonumber
\end{align}
where
\begin{equation}
\breve{v}_{\mathbf{n}}(t_{k})=(\breve{v}(t_{k},x))_{\mathbf{n}}=\left(
V(t_{k},x)\bar{v}(t_{k},x)\right)  _{\mathbf{n}}. \label{alg3}%
\end{equation}
To find $\breve{v}_{\mathbf{n}}(t_{k})$ one can either multiply two partial
sums of the form (\ref{al1}) and (\ref{ext}) or exploit fast Fourier transform
in the usual fashion (see, e.g. \cite{Can98}) to speed up the algorithm. The
algorithm (\ref{alg2}) can be viewed as analogous to spectral methods. It is
interesting that the layer method (\ref{NSM21})-(\ref{NSM23}) is, on the one
hand, related to a finite difference scheme (see below) and on the other hand,
to spectral methods.

\label{remFD}Let us discuss a relationship between the layer method
$(\ref{NSM21})$-$(\ref{NSM23})$ and finite difference methods. For simplicity
in writing, we give this illustration in the two-dimensional case. It is not
difficult to notice that the two-dimensional analog of the layer approximation
$(\ref{NSM21})$ can be re-written as the following finite difference scheme
for the SNSE $(\ref{NS1})$-$(\ref{NS3})$:
\begin{align}
&  \frac{\bar{v}(t_{k+1},x)-\bar{v}(t_{k},x)}{h}\label{fd1}\\
&  =\frac{\bar{v}(t_{k},x^{1}+\sigma\sqrt{h},x^{2}+\sigma\sqrt{h})+\bar
{v}(t_{k},x^{1}-\sigma\sqrt{h},x^{2}+\sigma\sqrt{h})-4\bar{v}(t_{k}%
,x^{1},x^{2})}{4h}\nonumber\\
&  +\frac{\bar{v}(t_{k},x^{1}+\sigma\sqrt{h},x^{2}-\sigma\sqrt{h})+\bar
{v}(t_{k},x^{1}-\sigma\sqrt{h},x^{2}-\sigma\sqrt{h})}{4h}\nonumber\\
&  -\frac{1}{\sigma\sqrt{h}}P\breve{v}(t_{k},x)+Pf(t_{k},x)+\sum_{r=1}%
^{q}\gamma_{r}(t_{k},x)\frac{\Delta w_{r}(t_{k+1})}{h}\nonumber
\end{align}
with
\begin{align}
\frac{\breve{v}(t_{k},x)}{\sigma\sqrt{h}}  &  =\bar{v}^{1}(t_{k},x)\frac
{\bar{v}(t_{k},x^{1}+\sigma\sqrt{h},x^{2}+\sigma\sqrt{h})-\bar{v}(t_{k}%
,x^{1}-\sigma\sqrt{h},x^{2}+\sigma\sqrt{h})}{4\sigma\sqrt{h}}\label{fd2}\\
&  +\bar{v}^{1}(t_{k},x)\frac{\bar{v}(t_{k},x^{1}+\sigma\sqrt{h},x^{2}%
-\sigma\sqrt{h})-\bar{v}(t_{k},x^{1}-\sigma\sqrt{h},x^{2}-\sigma\sqrt{h}%
)}{4\sigma\sqrt{h}}\nonumber\\
&  +\bar{v}^{2}(t_{k},x)\frac{\bar{v}(t_{k},x^{1}+\sigma\sqrt{h},x^{2}%
+\sigma\sqrt{h})-\bar{v}(t_{k},x^{1}+\sigma\sqrt{h},x^{2}-\sigma\sqrt{h}%
)}{4\sigma\sqrt{h}}\nonumber\\
&  +\bar{v}^{2}(t_{k},x)\frac{\bar{v}(t_{k},x^{1}-\sigma\sqrt{h},x^{2}%
+\sigma\sqrt{h})-\bar{v}(t_{k},x^{1}-\sigma\sqrt{h},x^{2}-\sigma\sqrt{h}%
)}{4\sigma\sqrt{h}}\ \ .\nonumber
\end{align}
As one can see, $\bar{v}(t_{k},\cdot)$ in the right-hand side of $(\ref{fd1})$
is evaluated at the nodes $(x^{1},x^{2}),$ $(x^{1}\pm\sigma\sqrt{h},x^{2}%
\pm\sigma\sqrt{h})$, which is typical for a standard explicit finite
difference scheme with the space discretization step $h_{x}$ taken equal to
$\sigma\sqrt{h}$ and $h$ being the time-discretization step. We also note that
if in the approximation $(\ref{NS30})$ we choose a different random vector
$\xi$ than in $(\ref{NS15})$ then we can obtain another layer method for the
SNSE which can be again re-written as a finite difference scheme (see such a
discussion in the case of the deterministic NSE in \cite{NS5}).

We recall \cite{M1,MT1,spde} that convergence theorems for layer methods (in
comparison with the theory of finite difference methods) do not contain any
conditions on stability of their approximations. In layer methods we do not
need to a priori prescribe space nodes: they are obtained automatically
depending on choice of a probabilistic representation and a numerical scheme.
We note that our error analysis for the layer methods (see Section~\ref{secER}%
) immediately implies the same error estimates for the corresponding finite
difference scheme $(\ref{fd1})$.

\begin{remark}
It is not difficult to see from $(\ref{fd2})$ that
\begin{equation}
(\bar{v}(t_{k},x),\nabla)\bar{v}(t_{k},x)\approx\frac{\breve{v}(t_{k}%
,x)}{\sigma\sqrt{h}}\ . \label{fd3}%
\end{equation}
If we put the exact $v(t_{k},x)$ in $(\ref{fd3})$ (both in its left and
right-hand sides) instead of the approximate $\bar{v}(t_{k},x)$ then the
accuracy of the approximation in $(\ref{fd3})$ is of order $O(h).$ This
observation is helpful for understanding a relationship between the layer
methods from this and the next section (see Remark~\ref{Remfd2} at the end of
the next section).
\end{remark}

\subsection{A layer method based on the direct probabilistic
representation\label{secDiL}}

The local version of probabilistic representation (\ref{NS11n})-(\ref{NS12n}) for the solution to the
backward SNSE (\ref{D1})-(\ref{D01}) has the form:
\begin{gather}
u(t_{k},x)=E^{\tilde{w}}u(t_{k+1},X_{t_{k},x}(t_{k+1}))\label{NS13n}\\
-E^{\tilde{w}}\left[  \int_{t_{k}}^{t_{k+1}}\{\nabla\tilde{p}(s^{\prime
},X_{t_{k},x}(s^{\prime}))-\tilde{f}(s^{\prime},X_{t_{k},x}(s^{\prime
})) \right. \nonumber \\
\Bigg.+(u(s^{\prime},X_{t_{k},x}(s^{\prime})),\nabla)u(s^{\prime},X_{t_{k}%
,x}(s^{\prime}))\}ds^{\prime}\Bigg] \nonumber\\
+\sum_{r=1}^{q}E^{\tilde{w}}\left[  \int_{t_{k}}^{t_{k+1}}\tilde{\gamma}%
_{r}(s^{\prime},X_{t_{k},x}(s^{\prime}))d\tilde{w}_{r}(s^{\prime})\right]
.\nonumber
\end{gather}

Using (\ref{NS13n}), we construct the one-step approximation of the solution to the
backward SNSE (\ref{D1})-(\ref{D01}):%
\begin{align}
u(t_{k},x)  &  =E^{\tilde{w}}u(t_{k+1},X_{t_{k},x}(t_{k+1}))-h\{\nabla
\tilde{p}(t_{k+1},x)-\tilde{f}(t_{k+1},x)\label{DL2} \displaybreak[0] \\
&  +(u(t_{k+1},x),\nabla)u(t_{k+1},x)\}+\sum_{r=1}^{q}\tilde{\gamma}%
_{r}(t_{k+1},x)\Delta_{k}\tilde{w}_{r}+\rho\nonumber \displaybreak[0]\\
&  =2^{-n}\sum_{j=1}^{2^{n}}u(t_{k+1},x+\sigma\sqrt{h}\xi_{j})\nonumber \displaybreak[0]\\
&  -h\{\nabla\tilde{p}(t_{k+1},x)-\tilde{f}(t_{k+1},x)+(u(t_{k+1}%
,x),\nabla)u(t_{k+1},x)\}\nonumber \displaybreak[0]\\
&  +\sum_{r=1}^{q}\tilde{\gamma}_{r}(t_{k+1},x)\Delta_{k}\tilde{w}_{r}%
+\rho,\nonumber
\end{align}
where $\rho=\rho(t_{k},x)$ is a remainder.

Using the Helmholtz-Hodge-Leray decomposition and taking into account that
$\operatorname{div}u(t_{k+1},x+\sigma\sqrt{h}\xi_{q})=0$ and
$\operatorname{div}\gamma_{r}=0,$\ we get from (\ref{DL2}):
\begin{gather}
u(t_{k},x)=2^{-n}\sum_{j=1}^{2^{n}}u(t_{k+1},x+\sigma\sqrt{h}\xi
_{j})-P[(u(t_{k+1},x),\nabla)u(t_{k+1},x)]h
\label{DL30}\\
+P\tilde{f}(t_{k+1},x)h-P^{\bot}[(u(t_{k+1},x),\nabla)u(t_{k+1},x)]h+P^{\bot}\tilde{f}(t_{k+1}%
,x)h\nonumber\\
-\nabla\tilde{p}(t_{k+1},x)h+\sum_{r=1}^{q}\tilde{\gamma}_{r}(t_{k+1},x)\Delta_{k}\tilde{w}_{r}%
+\rho,\nonumber
\end{gather}
whence we obtain after applying the operator $P:$%
\begin{align}
u(t_{k},x)  &  =2^{-n}\sum_{j=1}^{2^{n}}u(t_{k+1},x+\sigma\sqrt{h}\xi
_{j})-P[(u(t_{k+1},x),\nabla)u(t_{k+1},x)]h
\label{DL3}\\
&  +P\tilde{f}(t_{k+1},x)h+\sum_{r=1}^{q}\tilde{\gamma}_{r}(t_{k+1},x)\Delta_{k}\tilde{w}_{r}%
+P\rho.\nonumber
\end{align}
We re-write (\ref{DL30})-(\ref{DL3}) for the forward-time SNSE (\ref{NS1}%
)-(\ref{NS3}):
\begin{gather}
v(t_{k+1},x)=2^{-n}\sum_{j=1}^{2^{n}}v(t_{k},x+\sigma\sqrt{h}\xi
_{j})-P[(v(t_{k},x),\nabla)v(t_{k},x)]h\label{DL30n}\\
+Pf(t_{k},x)h-P^{\bot}[(v(t_{k},x),\nabla)v(t_{k},x)]h+P^{\bot}f(t_{k},x)h\nonumber\\
-\nabla
p(t_{k},x)h+\sum_{r=1}^{q}\gamma_{r}(t_{k},x)\Delta_{k}w_{r}+\rho\nonumber
\end{gather}
and
\begin{align}
v(t_{k+1},x)  &  =2^{-n}\sum_{j=1}^{2^{n}}v(t_{k},x+\sigma\sqrt{h}\xi
_{j})-P\left[  (v(t_{k},x),\nabla)v(t_{k},x)\right]  h\label{DL3n}\\
&  +Pf(t_{k},x)h+\sum_{r=1}^{q}\gamma_{r}(t_{k},x)\Delta_{k}w_{r}%
+P\rho.\nonumber
\end{align}
Dropping the remainder in (\ref{DL3n}), we obtain the one-step approximation
for the velocity $v(t_{k+1},x)$ in (\ref{NS1})-(\ref{NS3}):
\begin{align}
\hat{v}(t_{k+1},x)  &  =2^{-n}\sum_{j=1}^{2^{n}}v(t_{k},x+\sigma\sqrt{h}%
\xi_{j})-P\left[  (v(t_{k},x),\nabla)v(t_{k},x)\right]  h\label{DL5}\\
&  +Pf(t_{k},x)h+\sum_{r=1}^{q}\gamma_{r}(t_{k},x)\Delta_{k}w_{r}.\nonumber
\end{align}
It is easy to see that under Assumptions~2.1 $\operatorname{div}\hat
{v}(t_{k+1},x)=0.$ The corresponding layer method for the velocity of the SNSE
(\ref{NS1})-(\ref{NS3}) has the form%
\begin{gather}
\bar{v}(0,x)=\varphi(x), \label{DL7} \displaybreak[0]\\
\bar{v}(t_{k+1},x)=2^{-n}\sum_{j=1}^{2^{n}}\bar
{v}(t_{k},x+\sigma\sqrt{h}\xi_{j})-P\left[  (\bar{v}(t_{k},x),\nabla)\bar
{v}(t_{k},x)\right]  h \nonumber \displaybreak[0]\\
+Pf(t_{k},x)h+\sum_{r=1}^{q}\gamma_{r}(t_{k},x)\Delta_{k}w_{r},\ \ k=0,\ldots
,N-1.\nonumber
\end{gather}
This method can be turned into a numerical algorithm analogously to how we
constructed the numerical algorithm (\ref{alg2}) based on the layer method
(\ref{NSM21}) in Section~\ref{secSimL}.

\begin{remark}
\label{Remfd2}It is interesting to note (see also $(\ref{fd2})$ and
$(\ref{fd3})$) the relationship between the methods $(\ref{NSM21})$ and
$(\ref{DL7})$: $\sqrt{h}\breve{v}(t_{k},x)/\sigma$ from $(\ref{NSA4}%
)$-$(\ref{NSA6})$ is a finite-difference approximation of the term $(\bar
{v}(t_{k},x),\nabla)\bar{v}(t_{k},x)h$ in $(\ref{DL7})$. We remark that this
finite difference naturally arises via the probabilistic approach. It is
useful to have both methods in the arsenal of layer methods for SNSE: while
the method $(\ref{DL7})$ has a smaller one-step error than $(\ref{NSM21}),$ it
requires evaluation of spatial derivatives of $\bar{v}(t_{k},x)$.
\end{remark}

\subsection{Approximation of pressure\label{secPres}}

In the previous sections we constructed numerical methods for velocity
$v(t,x),$ in this section we propose approximations for pressure $p(t,x)$.

Applying the projection operator $P^{\bot}$ to SNSE (\ref{NS1})-(\ref{NS3}),
we get (see also (\ref{NS03})):
\begin{equation}
\nabla p(t,x)=-P^{\bot}\left[  (v(t,x),\nabla)v(t,x)\right]  +P^{\bot}f(t,x).
\label{pre1}%
\end{equation}
Based on (\ref{pre1}), we complement the layer method (\ref{DL7}) for the
velocity by the approximation of pressure as follows:
\begin{equation}
\nabla\bar{p}(t_{k+1},x)=-P^{\bot}\left[  (\bar{v}(t_{k+1},x),\nabla)\bar
{v}(t_{k+1},x)\right]  +P^{\bot}f(t_{k+1},x). \label{DL7p}%
\end{equation}

As a result, we obtain \textit{the layer method} (\ref{DL7}), (\ref{DL7p}) for
the solution of SNSE (\ref{NS1})-(\ref{NS3}).

It is clear that the numerical error $\nabla\bar{p}(t_{k+1},x)-\nabla p(t,x)$
is of the same order as the global errors of $\bar{v}(t_{k+1},x)$ and
$\nabla\bar{v}(t_{k+1},x).$ We note that in (\ref{DL7p}) to evaluate pressure
at time $t_{k+1}$ we use velocity at time $t_{k+1},$ i.e., the updated velocity.

\begin{remark}
We observe that $\rho$ in $(\ref{DL30n})$ is such that $P^{\bot}\rho=0.$
Indeed, it follows from $(\ref{DL30n})$-$(\ref{DL3n})$ (with $t_{k+1}$ instead
of $t_{k}$) that
\begin{equation}
\nabla p(t_{k+1},x)=-P^{\bot}\left[  (v(t_{k+1},x),\nabla)v(t_{k+1},x)\right]
+P^{\bot}f(t_{k+1},x)+P^{\bot}\rho. \label{pre2}%
\end{equation}
Comparing $(\ref{pre1})$ and $(\ref{pre2})$, we get $P^{\bot}\rho=0.$
\end{remark}

Let us now return to the layer method (\ref{NSM21}) for velocity. We have to
complement it with an approximation of pressure. To this end, we approximate
(see Remark~\ref{Remfd2} and (\ref{fd3})) the term $(\bar{v}(t_{k+1}%
,x),\nabla)\allowbreak\bar{v}(t_{k+1},x)$ in (\ref{DL7p}) by $\breve
{v}(t_{k+1},x)/\sigma\sqrt{h}$ with $\breve{v}(t_{k+1},x)$ from (\ref{NSM23})
(with $t_{k+1}$ instead of $t_{k})$. We obtain
\begin{equation}
\nabla\bar{p}(t_{k+1},x)=-\frac{1}{\sigma\sqrt{h}}P^{\bot}\breve{v}%
(t_{k+1},x)+P^{\bot}f(t_{k+1},x), \label{NSMp}%
\end{equation}
where $\breve{v}(t_{k+1},x)$ is from (\ref{NSM23}). Note that in the velocity
approximation (\ref{NSM21}) we use $\breve{v}(t_{k},x)$ while in the pressure
approximation (\ref{NSMp}) we use $\breve{v}(t_{k+1},x).$

As a result, we obtain \textit{the layer method} (\ref{NSM21})-(\ref{NSM23}),
(\ref{NSMp}) for the solution of SNSE (\ref{NS1})-(\ref{NS3}).

We remark that the layer method (\ref{NS18}) for velocity can be completed by
approximating the pressure as it is either in (\ref{DL7p}) with $\bar
{v}(t_{k+1},x)$ found due to (\ref{NS18}) or in (\ref{NSMp}) but with
$\breve{v}(t_{k+1},x)$ from (\ref{NSM23}) using $\bar{v}(t_{k+1},x)$ found
due to (\ref{NS18}).

To provide an example of an algorithm involving an approximation of pressure,
let us return to the algorithm (\ref{alg2}) for velocity. Based on
(\ref{NSMp}) (see also (\ref{N01})), we obtain
\begin{equation}
\bar{p}_{\mathbf{n}}(t_{k+1})=i\frac{L}{2\pi}\left(  \frac{\breve
{v}_{\mathbf{n}}^{\top}(t_{k+1})\mathbf{n}}{\sigma\sqrt{h}|\mathbf{n}|^{2}%
}-\frac{f_{\mathbf{n}}^{\top}(t_{k+1})\mathbf{n}}{|\mathbf{n}|^{2}}\right)
,\ \ \mathbf{n\neq0,\ }\bar{p}_{\mathbf{0}}(t_{k+1})=0, \label{algp}%
\end{equation}
where $\breve{v}_{\mathbf{n}}^{\top}(t_{k+1})$ are as in (\ref{alg3}) with
$t_{k+1}$ instead of $t_{k}$.

As a result, we obtain \textit{the algorithm} (\ref{alg2})-(\ref{alg3}),
(\ref{algp}) for the solution of SNSE (\ref{NS1})-(\ref{NS3}) which
corresponds to the layer method (\ref{NSM21})-(\ref{NSM23}), (\ref{NSMp}).

Analogously, one can obtain algorithms corresponding to the other two layer
methods considered in the paper.

\section{Error analysis\label{secER}}

In this section we provide theoretical support for the numerical methods from
the previous section. For definiteness, we consider the layer method
(\ref{NS18}). Analogous results can be obtained for the other two layer
methods proposed in Sections~\ref{secSimL}\ and~\ref{secDiL}.

As before, $||u(\cdot)||=||u(x)||$ denotes the $\mathbf{L}^{2}$-norm of a
function $u(x),$ $x\in Q.$ In this section we use the same letter $K$ for
various deterministic constants and $C=C(\omega)$ for various positive random variables.

We start with analysis of the local mean-square error.

\begin{theorem}
\label{lemonest}Let Assumptions~2.1 hold with $m_{0}>6$. The one-step error
\begin{equation}
\rho(t_{k+1},x)=\hat{v}(t_{k+1},x)-v(t_{k+1},x) \label{onesterr}%
\end{equation}
of the one-step approximation $(\ref{NSA1})$-$(\ref{NSA3})$ for the SNSE
$(\ref{NS1})$-$(\ref{NS3})$ is estimated as
\begin{equation}
||E(\rho(t_{k+1},x)|\mathcal{F}_{t_{k}}^{w})||\leq C(\omega)h^{2}%
,\ \label{lm1}%
\end{equation}
and for $1\leq p<p_{0}$
\begin{equation}
\left(  E||\rho(t_{k+1},\cdot)||^{2p}\right)  ^{1/2p}\leq Kh^{3/2}%
,\ \label{lm22}%
\end{equation}
where a random constant $C(\omega)>0$ with $EC^{2}<\infty$ does not depend on
$h$ and $k,$ a deterministic constant $K>0$ does not depend on $h$ and $k$ but
depends on $p,$ and $p_{0}=p_{0}(m_{0})>1$ is a positive number or
$p_{0}=\infty.$
\end{theorem}

\noindent\textbf{Proof. } Using Assumptions~2.1, we expand the right-hand side
of (\ref{NSA3}), substitute the outcome in (\ref{NSA1}), and obtain
\begin{align}
\hat{v}(t_{k+1},x)  &  =v(t_{k},x)-hP\left[  (v(t_{k},x),\nabla)v(t_{k}%
,x)\right]  +\frac{\sigma^{2}}{2}h\Delta v(t_{k},x)\label{lm3}\\
&  +Pf(t_{k},x)h+\sum_{r=1}^{q}\gamma_{r}(t_{k},x)\Delta_{k}w_{r}+r_{1}%
(t_{k},x),\nonumber
\end{align}
where the remainder $r_{1}(t_{k},x)$ has the form
\begin{align*}
r_{1}(t_{k},x)  &  =\frac{h^{2}}{2}\sum_{i,j=1}^{n}P\left[  v^{i}%
(t_{k},x)v^{j}(t_{k},x)\frac{\partial^{2}}{\partial x^{i}\partial x^{j}%
}v(t_{k},\Theta)\right] \\
&  +\frac{\sigma^{2}h^{2}}{2}\sum_{i,j=1}^{n}P\left[  v^{j}(t_{k}%
,x)\frac{\partial^{2}}{\left(  \partial x^{i}\right)  ^{2}\partial x^{j}%
}v(t_{k},\tilde{\Theta})\right] \\
&  +\frac{\sigma^{4}h^{2}}{24}2^{-n}\sum_{j=1}^{2^{n}}\sum_{i=1}^{n}P\left[
\frac{\partial^{4}}{\partial x^{i_{1}}\partial x^{i_{2}}\partial x^{i_{3}%
}\partial x^{i_{4}}}v(t_{k},\Xi_{j})\xi_{j}^{i_{1}}\xi_{j}^{i_{2}}\xi
_{j}^{i_{3}}\xi_{j}^{i_{4}}\right]  ,
\end{align*}
and $\Theta$ and $\tilde{\Theta}$ are some intermediate points between $x$ and
$x-v(t_{k},x)h,$ and $\Xi_{j}$ are some intermediate points between
$x-v(t_{k},x)h$ and $x-v(t_{k},x)h+\sigma\sqrt{h}\xi_{j}$ (we note that
$r_{1}$ is a vector and the intermediate points depend on the component of
$r_{1}$ but we do not reflect this in the notation). It is not difficult to
estimate that this remainder satisfies the inequalities
\begin{equation}
||E\left(  r_{1}(t_{k},x)|\mathcal{F}_{t_{k}}^{w}\right)  ||\leq
C(\omega)h^{2},\ \ \left(  E||r_{1}(t_{k},\cdot)||^{2p}\right)  ^{1/2p}\leq
Kh^{2}. \label{lm4}%
\end{equation}

We write the solution $v(s,x),\ s\geq t_{k},$ of (\ref{NS1})-(\ref{NS3}) as
\begin{align}
v(s,x)  &  =v(t_{k},x)+\int_{t_{k}}^{s}\left[  \frac{\sigma^{2}}{2}\Delta
v(s^{\prime},x)-(v(s^{\prime},x),\nabla)v(s^{\prime},x)+f(s^{\prime
},x)\right]  ds^{\prime}\label{lm5}\\
&  -\int_{t_{k}}^{s}\nabla p(s^{\prime},x)ds^{\prime}+\sum_{r=1}^{q}%
\int_{t_{k}}^{s}\gamma_{r}(s^{\prime},x)dw_{r}(s^{\prime})\nonumber
\end{align}
and, in particular,%
\begin{align}
v(t_{k+1},x)  &  =v(t_{k},x)+\int_{t_{k}}^{t_{k+1}}\left[  \frac{\sigma^{2}%
}{2}\Delta v(s,x)-(v(s,x),\nabla)v(s,x)+f(s,x)\right]  ds\label{lm6}\\
&  -\int_{t_{k}}^{t_{k+1}}\nabla p(s,x)ds+\sum_{r=1}^{q}\int_{t_{k}}^{t_{k+1}%
}\gamma_{r}(s,x)dw_{r}(s).\nonumber
\end{align}
Substituting $v(s,x)$ from (\ref{lm5}) in the integrand of the first integral
in (\ref{lm6}) and expanding $\gamma_{r}(s,x)$ at $(t_{k},x),$ we obtain
\begin{align}
v(t_{k+1},x)  &  =v(t_{k},x)+h\frac{\sigma^{2}}{2}\Delta v(t_{k}%
,x)-h(v(t_{k},x),\nabla)v(t_{k},x)+hf(t_{k},x)\label{lm7}\\
&  -\int_{t_{k}}^{t_{k+1}}\nabla p(s,x)ds+\sum_{r=1}^{q}\gamma_{r}%
(t_{k},x)\Delta_{k}w_{r}+r_{2}(t_{k},x),\nonumber
\end{align}
where
\[
r_{2}(t_{k},x)=r_{2}^{(1)}(t_{k},x)+r_{2}^{(2)}(t_{k},x)
\]
and
\begin{align*}
r_{2}^{(1)}(t_{k},x)   & =\frac{\sigma^{2}}{2}\int_{t_{k}}^{t_{k+1}}\left[
\int_{t_{k}}^{s}\Delta\left(  \frac{\sigma^{2}}{2}\Delta v(s^{\prime
},x)-(v(s^{\prime},x),\nabla)v(s^{\prime},x) \right. \right . \displaybreak[0]\\
& \bigg. \bigg. +f(s^{\prime},x)\bigg)ds^{\prime}\bigg]  ds
-\frac{\sigma^{2}}{2}\int_{t_{k}}^{t_{k+1}}\int_{t_{k}}^{s}\Delta\nabla
p(s^{\prime},x)ds^{\prime}ds \displaybreak[0]\\
&-\int_{t_{k}}^{t_{k+1}}(v(s,x),\nabla)
\left[  \int_{t_{k}}^{s}\left(
\frac{\sigma^{2}}{2}\Delta v(s^{\prime},x)
-(v(s^{\prime},x),\nabla)v(s^{\prime},x)
 \right. \right. \displaybreak[0]\\
&\bigg. \bigg. +f(s^{\prime},x)\bigg)  ds^{\prime}\bigg]  ds \displaybreak[0]\\
&  +\int_{t_{k}}^{t_{k+1}}(v(s,x),\nabla)\int_{t_{k}}^{s}\nabla p(s^{\prime
},x)ds^{\prime}ds \displaybreak[0]\\
&  -\int_{t_{k}}^{t_{k+1}}\left(  \int_{t_{k}}^{s}\left(  \frac{\sigma^{2}}%
{2}\Delta v(s^{\prime},x)-(v(s^{\prime},x),\nabla)v(s^{\prime},x) \right. \right. \displaybreak[0]\\
& \bigg. \bigg. +f(s^{\prime},x)\bigg)  ds^{\prime},\nabla\bigg) v(s,x)ds \displaybreak[0]\\
&+\int_{t_{k}}^{t_{k+1}}\left(  \int_{t_{k}}^{s}\nabla p(s^{\prime
},x)ds^{\prime},\nabla\right)  v(s,x)ds \displaybreak[0]\\
&  +\int_{t_{k}}^{t_{k+1}}(t_{k+1}-s)\frac{\partial}{\partial s}f(s,x)ds,
\end{align*}%
\begin{align*}
r_{2}^{(2)}(t_{k},x)  &  =\frac{\sigma^{2}}{2}\sum_{r=1}^{q}\int_{t_{k}%
}^{t_{k+1}}\int_{t_{k}}^{s}\Delta\gamma_{r}(s^{\prime},x)dw_{r}(s^{\prime
})ds \displaybreak[0]\\
&  -\sum_{r=1}^{q}\int_{t_{k}}^{t_{k+1}}\left[  (v(s,x),\nabla)\int_{t_{k}%
}^{s}\gamma_{r}(s^{\prime},x)dw_{r}(s^{\prime})\right]  ds \displaybreak[0]\\
&  -\sum_{r=1}^{q}\int_{t_{k}}^{t_{k+1}}\left(  \int_{t_{k}}^{s}\gamma
_{r}(s^{\prime},x)dw_{r}(s^{\prime}),\nabla\right)  v(s,x)ds \displaybreak[0]\\
&  +\sum_{r=1}^{q}\int_{t_{k}}^{t_{k+1}}\left(  w_{r}(t_{k+1})-w_{r}%
(s)\right)  \frac{\partial}{\partial s}\gamma_{r}(s,x)ds.
\end{align*}
We see that the remainder $r_{2}(t_{k},x)$ consists of 1) $r_{2}^{(1)}%
(t_{k},x)$ with terms of mean-square order $h^{2}$ and 2) $r_{2}^{(2)}%
(t_{k},x)$ with terms containing $\mathcal{F}_{t_{k+1}}^{w}$-measurable Ito
integrals of mean-square order $h^{3/2}$ which expectations with respect to
$\mathcal{F}_{t_{k}}^{w}$ equal zero.\ Further, using Assumptions~2.1, one can
show that
\begin{equation}
|E\left(  r_{2}(t_{k},x)|\mathcal{F}_{t_{k}}^{w}\right)  |\leq C(\omega
)h^{2},\ \ \left(  E\left\vert r_{2}(t_{k},x)\right\vert ^{2p}\right)
^{1/2p}\leq Kh^{3/2},\ \label{lm8}%
\end{equation}
where $C(\omega)>0$ and $K>0$ do not depend on $k,$ $x,$ and $h.$ Based on the
second inequality in (\ref{lm8}), we obtain
\begin{align}
E\left\vert |r_{2}(t_{k},\cdot)|\right\vert ^{2p}  &  =E\left(  \int%
_{Q}\left[  r_{2}(t_{k},x)\right]  ^{2}dx\right)  ^{p}\leq KE\int%
_{Q}\left\vert r_{2}(t_{k},x)\right\vert ^{2p}dx\label{lm82}\\
&  \leq K\int_{Q}E\left\vert r_{2}(t_{k},x)\right\vert ^{2p}dx\leq
Kh^{2p\times3/2}\ .\nonumber
\end{align}

Applying the projector operator $P$ to the left- and right-hand sides of
(\ref{lm7}), we arrive at
\begin{align}
v(t_{k+1},x)  &  =v(t_{k},x)+h\frac{\sigma^{2}}{2}\Delta v(t_{k}%
,x)-hP[(v(t_{k},x),\nabla)v(t_{k},x)]+hPf(t_{k},x)\label{lm9}\\
&  +\sum_{r=1}^{q}\gamma_{r}(t_{k},x)\Delta_{k}w_{r}+r_{3}(t_{k}%
,x),\ \nonumber
\end{align}
where the new remainder $r_{3}(t_{k},x)=Pr_{2}(t_{k},x).$ Using (\ref{lm82}),
we get
\begin{equation}
E||r_{3}(t_{k},\cdot)||^{2p}=E||Pr_{2}(t_{k},\cdot)||^{2p}\leq E||r_{2}%
(t_{k},\cdot)||^{2p}\leq Kh^{2p\times3/2}. \label{lm100}%
\end{equation}
Hence from here, (\ref{lm4}) and (\ref{lm3}), (\ref{lm9}), we obtain
(\ref{lm22}).

Observing that expectation of projection $P$ of Ito integrals remains equal to
zero, we get $E\left(  Pr_{2}^{(2)}(t_{k},x)|\mathcal{F}_{t_{k}}^{w}\right)
=0.$ Since $r_{2}^{(1)}(t_{k},x)$ consists of terms of mean-square order
$h^{2}$, we obtain%
\begin{align*}
||E\left(  r_{3}(t_{k},x)|\mathcal{F}_{t_{k}}^{w}\right)  ||^{2}  &
=||E\left(  Pr_{2}^{(1)}(t_{k},x)|\mathcal{F}_{t_{k}}^{w}\right)  ||^{2} \displaybreak[0]\\%
&=\int_{Q}\left[  E\left(  Pr_{2}^{(1)}(t_{k},x)|\mathcal{F}_{t_{k}}%
^{w}\right)  \right]  ^{2}dx \displaybreak[0]\\
&  \leq\int_{Q}E\left(  \left[  Pr_{2}^{(1)}(t_{k},x)\right]  ^{2}%
|\mathcal{F}_{t_{k}}^{w}\right)  dx \displaybreak[0]\\
&=E\left(  \int_{Q}\left[  Pr_{2}%
^{(1)}(t_{k},x)\right]  ^{2}dx|\mathcal{F}_{t_{k}}^{w}\right) \displaybreak[0]\\
&  \leq E\left(  \int_{Q}\left[  r_{2}^{(1)}(t_{k},x)\right]  ^{2}%
dx|\mathcal{F}_{t_{k}}^{w}\right)  \leq C(\omega)h^{4}%
\end{align*}
whence
\begin{equation}
||E\left(  r_{3}(t_{k},x)|\mathcal{F}_{t_{k}}^{w}\right)  ||\leq
C(\omega)h^{2}\ . \label{lm11}%
\end{equation}
Then the estimate (\ref{lm1}) follows from (\ref{lm4}), (\ref{lm11}) and
(\ref{lm3}), (\ref{lm9}). \ $\square$

\begin{remark}
We recall that in Assumptions~2.1 we require existence of moments of order
$m,$ $2\leq m<m_{0},$ of the solution and its spatial derivatives. The higher
the $m_{0},$ the higher $p$, $1\leq p<p_{0},$ can be taken in $(\ref{lm22})$.
In particular, to guarantee $(\ref{lm22})$ with $p=1,$ we need existence of
moments of up to the order $m=6,$ while if the moments of any order $m$ (i.e.,
$m_{0}=\infty)$ are finite then $(\ref{lm22})$ is valid for any $p.$ We also
note that the smoothness conditions on the SNSE solution (see Assumptions~2.1)
required for proving Theorem~\ref{lemonest} are so that $v(t,x)$ should have
continuous spatial derivatives up to order four and $p(t,x)$ -- up to order three.
\end{remark}

\begin{corollary}
\label{coras}Let Assumptions~2.1 hold with the bounded moments of any order
$m\geq2.$ Then for almost every trajectory $w(\cdot)$ and any $0<\varepsilon
<3/2$ there exists a constant $C(\omega)>0$ such that the one-step error from
$(\ref{onesterr})$ is estimated as
\begin{equation}
||\rho(t_{k+1},\cdot)||\leq C(\omega)h^{3/2-\varepsilon}, \label{coroe}%
\end{equation}
i.e., the layer method $(\ref{NS18})$ has the one-step error of order
$3/2-\varepsilon$ a.s.\ .
\end{corollary}

\noindent\textbf{Proof.} Here we follow the recipe used in
\cite{Gyo98,filter,spde}. The Markov inequality together with (\ref{lm22})
implies
\[
P(||\rho(t_{k+1},\cdot)||>h^{\gamma})\leq\frac{E||\rho(t_{k+1},\cdot)||^{2p}%
}{h^{2p\gamma}}\leq Kh^{2p(3/2-\gamma)}.
\]
Then for any $\gamma=3/2-\varepsilon$ there is a sufficiently large $p\geq1$
such that (recall that $h=T/N)$
\[
\sum_{N=1}^{\infty}P\left(  ||\rho(t_{k+1},\cdot)||>\frac{T^{\gamma}%
}{N^{\gamma}}\right)  \leq KT^{2p(3/2-\gamma)}\sum_{N=1}^{\infty}\frac
{1}{N^{2p(3/2-\gamma)}}<\infty.
\]
Hence, due to the Borel-Cantelli lemma, the random variable
$$\varsigma
:=\sup_{h>0}h^{-\gamma}||\rho(t_{k+1},\cdot)||$$
is a.s. finite which implies
(\ref{coroe}). \ $\square$

\begin{remark}
Since it is desirable for the order of the one-step error $||\rho
(t_{k+1},\cdot)||$ to be greater than one, we should impose the restriction on
$\varepsilon$ in $(\ref{coroe})$ to be in $(0,0.5).$ If we restrict ourselves
to fulfilment of the inequality $(\ref{coroe})$ with $\varepsilon
_{0}<\varepsilon<1/2,$ where $\varepsilon_{0}$ is some positive number, then
the conditions of Corollary~\ref{coras}\ can be weakened since for such
$\varepsilon$ it is sufficient to take $p_{0}=1/(2\varepsilon_{0}).\ $
\end{remark}

The intuition built on numerics for ordinary stochastic differential equations
(see, e.g. \cite{MT1}) and also based on layer methods for SPDEs
\cite{filter,spde} together with convergence results for layer methods for
deterministic NSE \cite{BM,NS5} suggests that the one-step error properties
proved in Theorem~\ref{lemonest} should lead to mean-square convergence of the
layer method (\ref{NS18}) with order one, i.e.,
\begin{equation}
(E||\bar{v}(t_{k},\cdot)-v(t_{k},\cdot)||^{2p})^{1/2p}\leq Kh. \label{msqone}%
\end{equation}
However, we have not succeeded in proving such a result. Below we prove an
almost sure (a.s.) convergence of the method (\ref{NS18}) with lower order of
$1/2-\varepsilon$ for arbitrary $\varepsilon>0$ than the $1-\varepsilon$ a.s.
order which should follow from (\ref{msqone}) and the Borel-Cantelli-type of
arguments (see, e.g. \cite{filter,spde} and also the proof of
Corollary~\ref{coras} above). In our numerical experiments (see
Section~\ref{secnum}) we observed the first order (both mean-square and a.s.)
convergence of a layer method on test examples.

Since we assumed in Assumptions~2.1\textit{ }that the problem\textit{
}(\ref{NS1})-(\ref{NS3})\textit{ }has a unique classical solution
$v(t,x),$\ $p(t,x)$ which has continuous derivatives in the space variable
$x$\ up to some order and since we are considering the periodic case, then
$v(t,x),\ p(t,x)$ and their derivatives are a.s. finite on $[0,T]\times Q$.

To prove the below a.s. convergence Theorem~\ref{tmhconuadd}, we make the
following assumptions on the approximate solution $\bar{v}(t_{k},x)$ from
(\ref{NS18}). \medskip

\noindent\textbf{Assumptions 4.1. }\textit{Let }$\bar{v}(t_{k},x),$
$k=0,\ldots,N,$ have continuous first-order spatial derivatives and
\begin{align}
|\bar{v}(t_{k},x)|  &  \leq C(\omega),\ \label{NS201}\\
|\partial\bar{v}(t_{k},x)/\partial x^{i}|  &  \leq C(\omega),\ \ i=1,\ldots
,n,\nonumber
\end{align}
\textit{where} $C(\omega)>0$\ \textit{is an a.s. finite constant independent
of} $x,\ h,\ k.$ \medskip

The first inequality in (\ref{NS201}) is necessary for a.s. convergence of the
layer method (\ref{NS18}). The second inequality is also necessary if one
expects convergence of spatial derivatives of $\bar{v}(t,x).$ We note that
even in the case of deterministic NSE \cite{BM,NS5} it turns out to be
problematic to derive the inequalities (\ref{NS201}) for the approximate
solutions. At the same time, verifying Assumptions~4.1 in numerical
experiments is straightforward. We also note that in the case of Oseen-Stokes
equations we succeeded in deriving such estimates for approximate solutions
and their spatial derivatives.

\begin{theorem}
\label{tmhconuadd}Let Assumptions~2.1 hold with the bounded moments of any
order $m\geq2$ and Assumptions~4.1 also hold. For almost every trajectory
$w(\cdot)$ and any $0<\varepsilon<1/2$ there exists a constant $C(\omega)>0$
such that
\begin{equation}
||\bar{v}(t_{k},\cdot)-v(t_{k},\cdot)||\leq C(\omega)h^{1/2-\varepsilon},
\label{thm1}%
\end{equation}
i.e., the layer method $(\ref{NS18})$ for the SNSE $(\ref{NS1})$-$(\ref{NS3})$
converges with order $1/2-\varepsilon$ a.s..
\end{theorem}

\noindent\textbf{Proof. } First, we note that it is easy to see that under
Assumptions~2.1 and ~4.1:
\begin{equation}
\operatorname{div}\bar{v}(t_{k},x)=0. \label{bvdiv}%
\end{equation}
Denote the error of the method (\ref{NS18})-(\ref{NS19}) on the $k$th layer
by
\[
\varepsilon(t_{k},x)=\bar{v}(t_{k},x)-v(t_{k},x).
\]
Due to (\ref{NS18}) and (\ref{NS19}), we obtain
\begin{align*}
\varepsilon(t_{k+1},x)+v(t_{k+1},x)  &  =\bar{v}(t_{k+1},x) \displaybreak[0]\\
&  =2^{-n}\sum_{j=1}^{2^{n}}P\bar{v}(t_{k},x-\bar{v}(t_{k},x)h+\sigma\sqrt
{h}\xi_{j})+Pf(t_{k},x)h \displaybreak[0]\\
&  +\sum_{r=1}^{q}\gamma_{r}(t_{k},x)\Delta_{k}w_{r}\\
&  =2^{-n}\sum_{j=1}^{2^{n}}Pv(t_{k},x-\bar{v}(t_{k},x)h+\sigma\sqrt{h}\xi
_{j}) \displaybreak[0]\\
&  +2^{-n}\sum_{j=1}^{2^{n}}P\varepsilon(t_{k},x-\bar{v}(t_{k},x)h+\sigma
\sqrt{h}\xi_{j}) \displaybreak[0]\\
&  +Pf(t_{k},x)h+\sum_{r=1}^{q}\gamma_{r}(t_{k},x)\Delta_{k}w_{r}.
\end{align*}
Using Assumptions~2.1, we obtain
\begin{equation}
v(t_{k},x-\bar{v}(t_{k},x)h+\sigma\sqrt{h}\xi_{j})=v(t_{k},x-v(t_{k}%
,x)h+\sigma\sqrt{h}\xi_{j})+r_{j}(t_{k},x), \label{thm2}%
\end{equation}
where
\begin{equation}
|r_{j}(t_{k},x)|\leq C(\omega)|\varepsilon(t_{k},x)|h \label{thm3}%
\end{equation}
and $C(\omega)$ is an a.s. finite random variable. Hence
\begin{align*}
\varepsilon(t_{k+1},x)+v(t_{k+1},x)    =2^{-n}\sum_{j=1}^{2^{n}}%
Pv(t_{k},x-v(t_{k},x)h+\sigma\sqrt{h}\xi_{j}) \displaybreak[0]\\
+2^{-n}\sum_{j=1}^{2^{n}} Pr_{j}(t_{k},x)
  +2^{-n}\sum_{j=1}^{2^{n}}P\varepsilon(t_{k},x-\bar{v}(t_{k},x)h+\sigma
\sqrt{h}\xi_{j}) \displaybreak[0]\\
  +Pf(t_{k},x)h+\sum_{r=1}^{q}\gamma_{r}(t_{k},x)\Delta_{k}w_{r}.
\end{align*}
Then we get
\begin{align}
\varepsilon(t_{k+1},x)&=2^{-n}\sum_{j=1}^{2^{n}}P\varepsilon(t_{k},x-\bar
{v}(t_{k},x)h+\sigma\sqrt{h}\xi_{j})+2^{-n}\sum_{j=1}^{2^{n}}Pr_{j}%
(t_{k},x) \label{thm5} \displaybreak[0]\\
&+\rho(t_{k+1},x), \nonumber
\end{align}
where $\rho(t_{k+1},x)$ is the error (see (\ref{onesterr})) of the one-step
approximation (\ref{NSA1})-(\ref{NSA3}) and this one-step error satisfies the
inequality (\ref{coroe}) from Corollary~\ref{coras}. It follows from
(\ref{thm5}), (\ref{thm3}) and (\ref{coroe}) that
\begin{align}
||\varepsilon(t_{k+1},\cdot)||    \leq & 2^{-n}\sum_{j=1}^{2^{n}}||P\varepsilon
(t_{k},\cdot-\bar{v}(t_{k},\cdot)h+\sigma\sqrt{h}\xi_{j})||+2^{-n}\sum
_{j=1}^{2^{n}}||Pr_{j}(t_{k},\cdot)||\label{thm7} \displaybreak[0]\\
&  +||\rho(t_{k+1},\cdot)||\nonumber \displaybreak[0]\\
 \leq &2^{-n}\sum_{j=1}^{2^{n}}||\varepsilon(t_{k},\cdot-\bar{v}(t_{k}%
,\cdot)h+\sigma\sqrt{h}\xi_{j})||+2^{-n}\sum_{j=1}^{2^{n}}||r_{j}(t_{k}%
,\cdot)|| \nonumber \displaybreak[0] \\
&+||\rho(t_{k+1},\cdot)||\nonumber \displaybreak[0]\\
\leq & 2^{-n}\sum_{j=1}^{2^{n}}||\varepsilon(t_{k},\cdot-\bar{v}(t_{k}%
,\cdot)h+\sigma\sqrt{h}\xi_{j})||+C(\omega)||\varepsilon(t_{k},\cdot
)||h \nonumber \displaybreak[0] \\
&+C(\omega)h^{3/2-\varepsilon}.\nonumber
\end{align}

Consider $\delta(x)=\varepsilon(t_{k},x-\bar{v}(t_{k},x)h+\sigma\sqrt{h}%
\xi_{j}).$ Due to Assumptions~4.1, the function $y(x)=x-\bar{v}(t_{k}%
,x)h+\sigma\sqrt{h}\xi_{j}$ is a differentiable function with continuous
partial derivatives. Furthermore, using Assumptions~4.1, one can show that for
sufficiently small $h>0$ the function $y(x)=x-\bar{v}(t_{k},x)h+\sigma\sqrt
{h}\xi_{j}$ is injective. Then, taking into account the $Q$-periodicity of
$\bar{v}(t_{k},x)$ and $\varepsilon^{i}(t_{k},x),$ we obtain
\begin{align*}
||\delta(\cdot)||^{2}  &  =\int_{Q}\sum_{i=1}^{n}\left[  \varepsilon^{i}%
(t_{k},x-\bar{v}(t_{k},x)h+\sigma\sqrt{h}\xi_{j})\right]  ^{2}dx\\
&  =\int_{Q}\sum_{i=1}^{n}\left[  \varepsilon^{i}(t_{k},y)\right]  ^{2}%
\frac{D(x^{1}\ldots x^{n})}{D(y^{1}\ldots y^{n})}dy.
\end{align*}
Due to Assumptions~4.1 and due to (\ref{bvdiv}), we get
\begin{align}
\frac{D(y^{1}\ldots y^{n})}{D(x^{1}\ldots x^{n})}  &  =\left\vert
\begin{array}
[c]{cccc}%
1-h\frac{\partial\bar{v}^{1}(t_{k},x)}{\partial x^{1}} & -h\frac{\partial
\bar{v}^{1}(t_{k},x)}{\partial x^{2}} & \cdots & -h\frac{\partial\bar{v}%
^{1}(t_{k},x)}{\partial x^{n}} \\
-h\frac{\partial\bar{v}^{2}(t_{k},x)}{\partial x^{1}} & 1-h\frac{\partial
\bar{v}^{2}(t_{k},x)}{\partial x^{2}} & \cdots & -h\frac{\partial\bar{v}%
^{2}(t_{k},x)}{\partial x^{n}} \\
\cdots & \cdots & \cdots & \cdots \\
-h\frac{\partial\bar{v}^{n}(t_{k},x)}{\partial x^{1}} & -h\frac{\partial
\bar{v}^{n}(t_{k},x)}{\partial x^{2}} & \cdots & 1-h\frac{\partial\bar{v}%
^{n}(t_{k},x)}{\partial x^{n}}%
\end{array}
\right\vert \label{thm8} \displaybreak[0]\\
&  =1+C(\omega)h^{2},\nonumber
\end{align}
where $C(\omega)$ is an a.s. finite random variable. Then, we also have%
\[
\dfrac{D(x^{1}\ldots x^{n})}{D(y^{1}\ldots y^{n})}=1+C(\omega)h^{2}.
\]

We obtain from (\ref{thm7}) and (\ref{thm8}):
\begin{equation}
||\varepsilon(t_{k+1},\cdot)||\leq||\varepsilon(t_{k},\cdot)||+C(\omega
)||\varepsilon(t_{k},\cdot)||h+C(\omega)h^{3/2-\varepsilon}, \label{thm9}%
\end{equation}
whence (\ref{thm1}) follows.$\ \ \square$

\begin{remark}
We recall that we have proved in Theorem~\ref{lemonest} that the mean and
mean-square one-step errors of the layer method $(\ref{NS18})$ (and
analogously of the other two layer methods from Section~\ref{secLayer}) are of
orders $O(h^{2})$ and $O(h^{3/2}),$ respectively. This has given us the basis
to argue that the methods from Section~\ref{secLayer} are of global
mean-square order one (see $(\ref{msqone})$). The same intuition implies that
if we incorporate terms of mean-square order $O(h^{3/2})$ and of mean order
$O(h^{2})$ in these first order methods (and thus make the mean-square
one-step errors to be of order $O(h^{2})$ and the mean errors of order
$O(h^{3}))$ then they become of global mean-square order $3/2.$ The required
Ito integrals of mean-square order $O(h^{3/2})$ can be simulated in the
constructive way (and hence these methods of order $3/2$ are constructive). In
the case of deterministic NSE (i.e., when $\gamma_{r}=0)$ such a method of
global mean-square order $3/2$ becomes of order two and coincides with the
corresponding layer method derived in \cite{NS5}.
\end{remark}

Let us now consider the error of the approximations of pressure considered in
Section~\ref{secPres}. In the next proposition we prove convergence of
pressure evaluated by (\ref{DL7p}),\ (\ref{NS18}). Analogously, one can prove
convergence of the other approximations of pressure derived in
Section~\ref{secPres}.

\begin{proposition}
\label{prp32}Let assumptions of Theorem~\ref{tmhconuadd} hold. In addition
assume that second-order spatial derivatives of the approximate solution are
a.s. finite: $|\partial^{2}\bar{v}(t_{k},x)/\partial x^{i}\partial
x^{j}|\allowbreak\leq C(\omega).$ Then for almost every trajectory $w(\cdot)$
and any $0<\varepsilon<1/3$ there exists a constant $C(\omega)>0$ such that
the approximate pressure $\bar{p}(t_{k},x)$\ from $(\ref{DL7p})$%
,\ $(\ref{NS18})$\ satisfies the following inequality%
\begin{equation}
\Vert\bar{p}(t_{k},\cdot)-p(t_{k},\cdot)\Vert\leq C(\omega)h^{1/3-\epsilon}.
\label{NS204}%
\end{equation}

\end{proposition}

\noindent\textbf{Proof}. We have
\begin{align}
\frac{\partial v^{i}}{\partial x^{j}}(t_{k},x)  &  =\frac{v^{i}(t_{k},x+\delta
e_{j})-v^{i}(t_{k},x-\delta e_{j})}{2\delta}+O(\delta^{2}),\label{NS25} \displaybreak[0]\\
\frac{\partial\bar{v}^{i}}{\partial x^{j}}(t_{k},x)  &  =\frac{\bar{v}%
^{i}(t_{k},x+\delta e_{j})-\bar{v}^{i}(t_{k},x-\delta e_{j})}{2\delta
}+O(\delta^{2}),\nonumber
\end{align}
where $\delta$ is a positive sufficiently small number and $|O(\delta
^{2})|\leq C(\omega)\delta^{2}$. Due to Theorem~\ref{tmhconuadd},
\begin{gather}
\left\Vert \frac{v(t_{k},x+\delta e_{j})-v(t_{k},x-\delta e_{j})}{2\delta
}-\frac{\bar{v}(t_{k},x+\delta e_{j})-\bar{v}(t_{k},x-\delta e_{j})}{2\delta
}\right\Vert \label{NS255} \\
\leq C(\omega)\frac{h^{1/2-\epsilon/2}}{\delta}\ \text{a.s.} \nonumber
\end{gather}
Choosing $\delta=ch^{1/6+\epsilon/2}$ with some $c>0,$ we obtain from
(\ref{NS25}) and (\ref{NS255}) that%
\begin{equation}
\left\Vert \frac{\partial v}{\partial x^{j}}(t_{k},\cdot)-\frac{\partial
\bar{v}}{\partial x^{j}}(t_{k},\cdot)\right\Vert \leq C(\omega)h^{1/3-\epsilon
}\ \ \ \text{a.s.\ .} \label{NS26}%
\end{equation}
Subtracting (\ref{pre1}) with $t=t_{k}$ from (\ref{DL7p}) with $t_{k}$ instead
of $t_{k+1}$, we get%
\begin{gather}
\left\Vert \nabla\bar{p}(t_{k},\cdot)-\nabla p(t_{k},\cdot)\right\Vert
=\left\Vert P^{\bot}\left[  (v(t_{k},\cdot),\nabla)v(t_{k},\cdot)\right]
-P^{\bot}\left[  (\bar{v}(t_{k},\cdot),\nabla)\bar{v}(t_{k},\cdot)\right]
\right\Vert \label{NS27} \displaybreak[0]\\
\leq\left\Vert P^{\bot}\left[  (v(t_{k},\cdot),\nabla)(v(t_{k},\cdot)-\bar
{v}(t_{k},\cdot))\right]  \right\Vert +\left\Vert P^{\bot}\left[
(v(t_{k},\cdot)-\bar{v}(t_{k},\cdot),\nabla)\bar{v}(t_{k},\cdot)\right]
\right\Vert \nonumber \displaybreak[0]\\
\leq\left\Vert (v(t_{k},\cdot),\nabla)(v(t_{k},\cdot)-\bar{v}(t_{k}%
,\cdot))\right\Vert +\left\Vert (v(t_{k},\cdot)-\bar{v}(t_{k},\cdot
),\nabla)\bar{v}(t_{k},\cdot)\right\Vert \ .\nonumber
\end{gather}
Due to Assumptions~2.1 and (\ref{NS26}),
\begin{equation}
\left\Vert (v(t_{k},\cdot),\nabla)(v(t_{k},\cdot)-\bar{v}(t_{k},\cdot
))\right\Vert \leq C(\omega)h^{1/3-\epsilon}\ \ \text{a.s.\ .} \label{NS28}%
\end{equation}
Due to Assumptions~4.1 and Theorem~\ref{tmhconuadd},
\begin{equation}
\left\Vert (v(t_{k},\cdot)-\bar{v}(t_{k},\cdot),\nabla)\bar{v}(t_{k}%
,\cdot)\right\Vert \leq C(\omega)h^{1/2-\epsilon}\ \ \text{a.s.\ .}
\label{NS29}%
\end{equation}
Thus, (\ref{NS27})-(\ref{NS29}) imply (\ref{NS204}). \ $\ \square$

\begin{remark}
To prove the estimate
\begin{equation}
\left\Vert \frac{\partial v}{\partial x^{j}}(t_{k},x)-\frac{\bar{v}%
(t_{k},x+\delta e_{j})-\bar{v}(t_{k},x-\delta e_{j})}{2\delta}\right\Vert \leq
C(\omega)h^{1/3-\epsilon}\ \ \ \text{a.s.\ ,} \label{eremp1}%
\end{equation}
we do not need in the assumption on boundedness of second-order spatial
derivatives of the approximate solution. Then, under the conditions of
Theorem~\ref{tmhconuadd} (without the additional assumption on second-order
spatial derivatives of the approximate solution), we can analogously prove
convergence with a.s. order $1/3-\epsilon$ of the approximate pressure
$\bar{p}(t_{k},x)$\ from $(\ref{NSMp})$ with $\breve{v}(t_{k+1},x)$ from
$(\ref{NSM23})$ in which we substitute $\bar{v}(t_{k+1},x)$ found due to
$(\ref{NS18})$.
\end{remark}

\begin{remark}
\label{remp2}As we discussed earlier in this section, though we proved
$1/2-\varepsilon$ a.s. convergence order for the velocity approximation in
Theorem~\ref{tmhconuadd}, we are expecting that the actual a.s. convergence
order is $1-\varepsilon$ which was observed in our numerical experiments in
Section~\ref{secnum}. Analogously, we expect that spatial derivatives of the
approximate velocity converge with a.s. order $1-\varepsilon$ instead of
$1/3-\epsilon$ shown in $(\ref{NS26})$. It is not difficult to see from the
proof of Proposition~\ref{prp32} that a.s. convergence of both velocity and
its first-order spatial derivatives with order $1-\varepsilon$ implies a.s.
convergence of pressure with order $1-\varepsilon.$ In our numerical
experiments (see Section~\ref{secnum}) we observed convergence (both
mean-square and a.s.) of pressure with order one.
\end{remark}

\section{Numerical examples\label{secnum}}

In this section we test the numerical algorithm (\ref{alg2}) from
Section~\ref{secSimL} on two model problems. The experiments indicate that the
algorithm has the first order mean-square convergence.

\subsection{Model problems}

We introduce two model examples of SNSE (\ref{NS1})-(\ref{NS3}) which
solutions can be written in an analytic form. Both examples are
generalizations of the deterministic model of laminar flow from \cite{Taylor}
to the stochastic case.\medskip

\textbf{First model problem. }Let
\begin{equation}
f(t,x)=0,\text{\ \ }\varphi(x)=0, \label{nl1}%
\end{equation}%
\begin{align}
q  &  =1,\label{nl11}  \displaybreak[0]\\
\gamma_{1}^{1}(t,x)  &  =A\sin\frac{2\pi\kappa\ x^{1}}{L}\cos\frac{2\pi
\kappa\ x^{2}}{L}\exp\left(  -\sigma^{2}\left(  \frac{2\pi\kappa}{L}\right)
^{2}t\right)  \ ,\nonumber \displaybreak[0]\\
\gamma_{1}^{2}(t,x)  &  =-A\cos\frac{2\pi\kappa\ x^{1}}{L}\sin\frac{2\pi
\kappa\ x^{2}}{L}\exp\left(  -\sigma^{2}\left(  \frac{2\pi\kappa}{L}\right)
^{2}t\right)  \ ,\ \kappa\in\mathbf{Z},\ \ A\in\mathbf{R},\nonumber
\end{align}
then it is easy to check that the problem (\ref{NS1})-(\ref{NS3}),
(\ref{nl1})-(\ref{nl11}) has the following solution
\begin{align}
v^{1}(t,x)    =&A\sin\frac{2\pi\kappa\ x^{1}}{L}\cos\frac{2\pi\kappa\ x^{2}%
}{L} \label{nl2} \displaybreak[0] \\
&\times \exp\left(  -\sigma^{2}\left(  \frac{2\pi\kappa}{L}\right)  ^{2}t\right)
w(t)\ , \nonumber \displaybreak[0] \\
v^{2}(t,x)   = &-A\cos\frac{2\pi\kappa\ x^{1}}{L}\sin\frac{2\pi\kappa\ x^{2}%
}{L} \nonumber \displaybreak[0] \\
&\times \exp\left(  -\sigma^{2}\left(  \frac{2\pi\kappa}{L}\right)  ^{2}t\right)
w(t)\ ,\nonumber \displaybreak[0]\\
p(t,x)    =&\frac{A^{2}}{4}\left(  \cos\frac{4\pi\kappa\ x^{1}}{L}+\cos
\frac{4\pi\kappa\ x^{2}}{L}\right)  \nonumber \displaybreak[0] \\
& \times \exp\left(  -2\sigma^{2}\left(  \frac
{2\pi\kappa}{L}\right)  ^{2}t\right)  (w(t))^{2}\ .\nonumber
\end{align}

\textbf{Second model problem. }To construct this example, we recall the
following proposition from \cite{HRoz07}.

\begin{proposition}
\label{PropRoz}Let $V(t,x),$ $P(t,x)$ be a solution of the deterministic NSE
with zero forcing $($i.e., of $(\ref{NS1})$-$(\ref{NS3})$ with all $\gamma
_{r}=0$ and $f(t,x)=0)$ then the solution $v(t,x),$ $p(t,x)$ of $(\ref{NS1}%
)$-$(\ref{NS3})$ with constant $\gamma_{r}(t,x)=\gamma_{r}$ and $f(t,x)=0$ is
equal to
\begin{align}
v(t,x)  &  =V\left(  t,x-\int_{0}^{t}\sum_{r=1}^{q}\gamma_{r}w_{r}%
(s)ds\right)  +\sum_{r=1}^{q}\gamma_{r}w_{r}(t),\label{sgt}\\
p(t,x)  &  =P\left(  t,x-\int_{0}^{t}\sum_{r=1}^{q}\gamma_{r}w_{r}%
(s)ds\right)  . \label{sgt1}%
\end{align}

\end{proposition}

Combining this proposition with the deterministic model of laminar flow from
\cite{Taylor}, we obtain that if
\begin{gather}
f(t,x)=0,\text{\ \ }\varphi(x)=\left(  A\sin\frac{2\pi\kappa\ x^{1}}{L}%
\cos\frac{2\pi\kappa\ x^{2}}{L},-A\cos\frac{2\pi\kappa\ x^{1}}{L}\sin
\frac{2\pi\kappa\ x^{2}}{L}\right)  ^{\top},\label{nl3} \displaybreak[0]\\
\kappa\in\mathbf{Z},\ \ A\in\mathbf{R},\nonumber
\end{gather}
and%
\begin{equation}
q=1,\ \gamma_{1}^{1}(t,x)=\gamma^{1},\ \ \gamma_{1}^{2}(t,x)=\gamma^{2}\ .
\label{nl31}%
\end{equation}
then the problem (\ref{NS1})-(\ref{NS3}), (\ref{nl3})-(\ref{nl31}) has the
following solution
\begin{align}
v^{1}(t,x)  &  =A\sin\frac{2\pi\kappa\ \left(  x^{1}-\gamma^{1}I(t)\right)
}{L}\cos\frac{2\pi\kappa\ \left(  x^{2}-\gamma^{2}I(t)\right)  }{L} \label{nl4} \displaybreak[0]\\
& \times \exp\left(
-\sigma^{2}\left(  \frac{2\pi\kappa}{L}\right)  ^{2}t\right)
+\gamma^{1}w(t),\nonumber \displaybreak[0]\\
v^{2}(t,x)  &  =-A\cos\frac{2\pi\kappa\ \left(  x^{1}-\gamma^{1}I(t)\right)
}{L}\sin\frac{2\pi\kappa\ \left(  x^{2}-\gamma^{2}I(t)\right)  }{L} \nonumber \displaybreak[0]\\
& \times \exp\left(
-\sigma^{2}\left(  \frac{2\pi\kappa}{L}\right)  ^{2}t\right)
+\gamma^{2}w(t),\nonumber \displaybreak[0]\\
p(t,x)  &  =\frac{A^{2}}{4}\left(  \cos\frac{4\pi\kappa\ \left(  x^{1}%
-\gamma^{1}I(t)\right)  }{L}+\cos\frac{4\pi\kappa\ \left(  x^{2}-\gamma
^{2}I(t)\right)  }{L}\right)  \nonumber \displaybreak[0] \\
& \times \exp\left(  -2\sigma^{2}\left(  \frac{2\pi
\kappa}{L}\right)  ^{2}t\right)  ,\nonumber
\end{align}
where
\[
I(t)=\int_{0}^{t}w(s)ds,\ w(s)=w_{1}(s).
\]

\subsection{Results of numerical experiments}

In our numerical experiments we test the algorithm (\ref{alg2})-(\ref{alg3}),
(\ref{algp}) which is a realization of the layer method (\ref{NSM21}%
)-(\ref{NSM23}), (\ref{NSMp}). This algorithm possesses the following properties.

\begin{proposition}
\label{Propcons}\textbf{1.} The approximate solution of the problem
$(\ref{NS1})$-$(\ref{NS3})$, $(\ref{nl1})$-$(\ref{nl11})$ obtained by the
algorithm $(\ref{alg2})$-$(\ref{alg3})$, $(\ref{algp})$ contains only those
modes which are present in the coefficient $\gamma_{1}(t,x)$ from
$(\ref{nl11})$, i.e., which are present in the exact solution $(\ref{nl2})$.

\textbf{2.} The approximate solution of the problem $(\ref{NS1})$%
-$(\ref{NS3})$, $(\ref{nl3})$-$(\ref{nl31})$ obtained by the algorithm
$(\ref{alg2})$-$(\ref{alg3})$, $(\ref{algp})$ contains only those modes which
are present in the initial condition $\varphi(x)$ from $(\ref{nl3})$ and the
zero mode, i.e., which are present in the exact solution $(\ref{nl4})$.
\end{proposition}

The proof of this proposition is analogous to the proof of a similar result in
the deterministic case \cite{NS5} and it is omitted here. \medskip

We measure the numerical error in the experiments as follows. First, we
consider the relative mean-square error defined as%
\begin{equation}
err_{msq}^{v}=\frac{\sqrt{E\sum_{\mathbf{n}}|\bar{v}_{\mathbf{n}%
}(T)-v_{\mathbf{n}}(T)|^{2}}}{\sqrt{E\sum_{\mathbf{n}}|v_{\mathbf{n}}(T)|^{2}%
}}\ ,\ \ \ err_{msq}^{p}=\frac{\sqrt{E\sum_{\mathbf{n}}|\bar{p}_{\mathbf{n}%
}(T)-p_{\mathbf{n}}(T)|^{2}}}{\sqrt{E\sum_{\mathbf{n}}|p_{\mathbf{n}}(T)|^{2}%
}}\ . \label{msqerr}%
\end{equation}
Analysis of this error provides us with information about mean-square
convergence of the numerical algorithm considered. To evaluate this error in
the experiments, we use the Monte Carlo technique for finding the expectations
in (\ref{msqerr}) by running $K$ independent (with respect to realizations of
the Wiener process $w(t))$ realizations of $\bar{v}_{\mathbf{n}}(T),$
$v_{\mathbf{n}}(T),$ $\bar{p}_{\mathbf{n}}(T),\ p_{\mathbf{n}}(T).$ Second, we
consider the relative $L_{2}$-error for a fixed trajectory of $w(t):$
\begin{equation}
err^{v}=\frac{\sqrt{\sum_{\mathbf{n}}|\bar{v}_{\mathbf{n}}(T)-v_{\mathbf{n}%
}(T)|^{2}}}{\sqrt{\sum_{\mathbf{n}}|v_{\mathbf{n}}(T)|^{2}}}\ ,\ \ \ err^{p}%
=\frac{\sqrt{\sum_{\mathbf{n}}|\bar{p}_{\mathbf{n}}(T)-p_{\mathbf{n}}(T)|^{2}%
}}{\sqrt{\sum_{\mathbf{n}}|p_{\mathbf{n}}(T)|^{2}}}\ . \label{err}%
\end{equation}
Analysis of this error provides us with information about a.s. convergence of
the numerical algorithm. To evaluate this error in the tests, we fix a
trajectory $w(t),$ $0\leq t\leq T,$ which is obtained with a small time step.

We note that in the case of the considered examples and the tested algorithm
(see Proposition~\ref{Propcons}) $v_{\mathbf{n}}(T)$ are nonzero only for
$|\mathbf{n}^{1}|=|\mathbf{n}^{2}|=|\kappa|$ and $p_{\mathbf{n}}(T)$ are
nonzero only for $|\mathbf{n}^{1}|=2|\kappa|,$\ $\mathbf{n}^{2}=0$ and
$\mathbf{n}^{1}=0,$\ $|\mathbf{n}^{2}|=2|\kappa|$. Hence, the sums in
(\ref{msqerr}) and (\ref{err}) are finite here. This also implies that it is
sufficient here to take the cut-off parameter $M$ in the algorithm
(\ref{alg2})-(\ref{alg3}), (\ref{algp}) to be equal to $2|\kappa|.$

The test results for the algorithm (\ref{alg2})-(\ref{alg3}), (\ref{algp})
applied to the first model problem (\ref{NS1})-(\ref{NS3}), (\ref{nl1}%
)-(\ref{nl11}) are presented in Tables$~$\ref{tab1} and~\ref{tab2}. In
Table$~$\ref{tab1} the \textquotedblleft$\pm$\textquotedblright\ reflects the
Monte Carlo errors in evaluating of $err_{msq}^{v}$ and $err_{msq}^{p}$, they
give the confidence intervals for the corresponding values with probability
$0.95$.%

%TCIMACRO{\TeXButton{B}{\begin{table}[htb] \centering}}%
%BeginExpansion
\begin{table}[htb] \centering
%EndExpansion
\caption{Mean-square relative errors $err_{msq}^v$ and $err_{msq}^p$ from (\ref{msqerr}) at $T=3$ in simulation of the problem
(\ref{NS1})-(\ref{NS3}),  (\ref{nl1})-(\ref{nl11})
with $\sigma =0.1$, $A=1$, $\kappa =1$, $L=1$
by the algorithm (\ref{alg2})-(\ref{alg3}), (\ref{algp})  with $M=2$ and various time steps $h$.
The \textquotedblleft $\pm $\textquotedblright   reflects
the Monte Carlo error in evaluating $err_{msq}^v$ and $err_{msq}^p$  via the Monte Carlo technique
with $K=4000$ independent runs.
The exact values (up to 5 d.p.) of the denominators in (\ref{msqerr}) are $0.37470$ and $0.12159$, respectively.
\label{tab1}} \setlength{\tabcolsep}{3pt}%
\begin{tabular}
[c]{lc}\hline
$h$ &
\begin{tabular}
[c]{ll}%
velocity & \ \ \ \ \ \ \ \ \ \ \ \ \ \ \ \ pressure
\end{tabular}
\\\hline
$0.2$ & \multicolumn{1}{l}{%
\begin{tabular}
[c]{ll}%
$0.0537$ $\pm\ 0.0012$ & \ \ \ \ $\ \ 0.0710$ $\pm\ 0.0038$%
\end{tabular}
}\\
$0.1$ & \multicolumn{1}{l}{%
\begin{tabular}
[c]{ll}%
$0.0263$ $\pm\ 0.0006$ & \ \ \ \ $\ \ 0.0337$ $\pm\ 0.0016$%
\end{tabular}
}\\
$0.05$ & \multicolumn{1}{l}{%
\begin{tabular}
[c]{ll}%
$0.0130$ $\pm\ 0.0003$ & \ \ \ \ $\ \ 0.0170$ $\pm\ 0.0009$%
\end{tabular}
}\\
$0.02$ & \multicolumn{1}{l}{%
\begin{tabular}
[c]{ll}%
$0.0052$ $\pm\ 0.0001$ & \ \ \ \ $\ \ 0.0066$ $\pm\ 0.0003$%
\end{tabular}
}\\
$0.01$ & \multicolumn{1}{l}{%
\begin{tabular}
[c]{ll}%
$0.0025\ \pm$ $0.00006$ & \ \ \ $\ 0.0031$ $\pm\ 0.0001$%
\end{tabular}
}\\\hline
\end{tabular}%
%TCIMACRO{\TeXButton{E}{\end{table}}}%
%BeginExpansion
\end{table}%
%EndExpansion
%

%TCIMACRO{\TeXButton{B}{\begin{table}[htb] \centering}}%
%BeginExpansion
\begin{table}[htb] \centering
%EndExpansion
\caption{Relative errors $err^v$ and $err^p$ from (\ref{err}) at $T=3$ in simulation of the problem
(\ref{NS1})-(\ref{NS3}),  (\ref{nl1})-(\ref{nl11})
with $\sigma =0.1$, $A=1$, $\kappa =1$, $L=1$ for a fixed trajectory of the Wiener process $w(t)$
by the algorithm (\ref{alg2})-(\ref{alg3}), (\ref{algp})  with $M=2$ and various time steps $h$.
The exact values (up to 5 d.p.) of the denominators in (\ref{err}) are  $0.43950$ and $0.09658$, respectively.
\label{tab2}} \setlength{\tabcolsep}{3pt}%
\begin{tabular}
[c]{lc}\hline
$h$ &
\begin{tabular}
[c]{ll}%
velocity & \ \ \ pressure
\end{tabular}
\\\hline
$0.2$ & \multicolumn{1}{l}{%
\begin{tabular}
[c]{ll}%
$0.0485$ & \ \ \ \ $\ \ 0.0585$%
\end{tabular}
}\\
$0.1$ & \multicolumn{1}{l}{%
\begin{tabular}
[c]{ll}%
$0.0237$ & \ \ \ \ $\ \ 0.0284$%
\end{tabular}
}\\
$0.05$ & \multicolumn{1}{l}{%
\begin{tabular}
[c]{ll}%
$0.0117$ & \ \ \ \ $\ \ 0.0141$%
\end{tabular}
}\\
$0.02$ & \multicolumn{1}{l}{%
\begin{tabular}
[c]{ll}%
$0.0047$ & \ \ \ \ $\ \ 0.0056$%
\end{tabular}
}\\
$0.01$ & \multicolumn{1}{l}{%
\begin{tabular}
[c]{ll}%
$0.0023$ & \ \ \ \ $\ \ 0.0028$%
\end{tabular}
}\\\hline
\end{tabular}%
%TCIMACRO{\TeXButton{E}{\end{table}}}%
%BeginExpansion
\end{table}%
%EndExpansion

We can conclude from Table$~$\ref{tab1} that both velocity and pressure found
due to the algorithm (\ref{alg2})-(\ref{alg3}), (\ref{algp}) demonstrate the
mean-square convergence with order $1.$ We also see from Table$~$\ref{tab2}
that both velocity and pressure converge with order $1$ for a particular,
fixed trajectory of $w(t).$ We note that we repeated the experiment for other
realizations of $w(t)$ and observed the same behavior. The observed first
order convergence of the algorithm is consistent with our prediction (see
(\ref{msqone}), the discussion after it, and Remark~\ref{remp2}).

The test results for the algorithm (\ref{alg2})-(\ref{alg3}), (\ref{algp})
applied to the second model problem (\ref{NS1})-(\ref{NS3}), (\ref{nl3}%
)-(\ref{nl31}) are presented in Table$~$\ref{tab3}. In these tests we limit
ourselves to simulation for a particular, fixed trajectory of $w(t)$ and
observation of a.s. convergence. We note that evaluation of the exact solution
(\ref{nl4}) requires simulation of the integral $I(t).$ This was done in the
following way. At each time step $k+1,$ $k=0,\ldots,N-1,$ we simulate a Wiener
increment $\Delta_{k}w$ as i.i.d. Gaussian $\mathcal{N}(0,h)$ random variables
(and we find $w(t_{k+1})=w(t_{k})+\Delta_{k}w)$ and i.i.d. Gaussian
$\mathcal{N}(0,1)$ random variables $\eta_{k}.$ Then (see \cite[Chapter
1]{MT1}):
\[
I(t_{k+1})=I(t_{k})+hw(t_{k})+\frac{h}{2}\Delta_{k}w+\frac{h^{3/2}}{\sqrt{12}%
}\eta_{k}\ .
\]
%

%TCIMACRO{\TeXButton{B}{\begin{table}[htb] \centering}}%
%BeginExpansion
\begin{table}[htb] \centering
%EndExpansion
\caption{Relative errors $err^v$ and $err^p$ from (\ref{err}) at $T=3$ in simulation of the problem
(\ref{NS1})-(\ref{NS3}),  (\ref{nl3})-(\ref{nl31})
with $\sigma =0.1$, $A=1$, $\kappa =1$, $L=1$, $\gamma^1 =0.5$, $\gamma^2 =0.2$
for a fixed trajectory of the Wiener process $w(t)$
by the algorithm (\ref{alg2})-(\ref{alg3}), (\ref{algp})  with $M=2$ and various time steps $h$.
The exact values (up to 6 d.p.) of the denominators in (\ref{err}) are  $0.505620$ and $0.000548$, respectively.
\label{tab3}} \setlength{\tabcolsep}{3pt}%
\begin{tabular}
[c]{lc}\hline
$h$ &
\begin{tabular}
[c]{ll}%
velocity & \ \ \ pressure
\end{tabular}
\\\hline
$0.01$ & \multicolumn{1}{l}{%
\begin{tabular}
[c]{ll}%
$0.166\ \ $ & $\ \ \ \ \ \ 0.973$%
\end{tabular}
}\\
$0.005$ & \multicolumn{1}{l}{%
\begin{tabular}
[c]{ll}%
$0.068\ \ $ & $\ \ \ \ \ \ 0.384$%
\end{tabular}
}\\
$0.002$ & \multicolumn{1}{l}{%
\begin{tabular}
[c]{ll}%
$0.024\ \ $ & $\ \ \ \ \ \ 0.134$%
\end{tabular}
}\\
$0.001$ & \multicolumn{1}{l}{%
\begin{tabular}
[c]{ll}%
$0.0118$ & $\ \ \ \ \ \ 0.0645$%
\end{tabular}
}\\
$0.0005$ & \multicolumn{1}{l}{%
\begin{tabular}
[c]{ll}%
$0.0058$ & $\ \ \ \ \ \ 0.0313$%
\end{tabular}
}\\\hline
\end{tabular}%
%TCIMACRO{\TeXButton{E}{\end{table}}}%
%BeginExpansion
\end{table}%
%EndExpansion

Again, the observed first order convergence of the algorithm in
Table~\ref{tab3} is consistent with our prediction (see the discussion after
(\ref{msqone}) and Remark~\ref{remp2}).

\section*{Acknowledgments}

The work was partially supported by the Royal Society International Joint
Project grant JP091142.

\end{document}